\documentclass[reqno]{amsart}
\makeatletter
\def\oversortoftilde#1{\mathop{\vbox{\m@th\ialign{##\crcr\noalign{\kern3\p@}%
      \sortoftildefill\crcr\noalign{\kern3\p@\nointerlineskip}%
      $\hfil\displaystyle{#1}\hfil$\crcr}}}\limits}

\def\sortoftildefill{$\m@th \setbox\z@\hbox{$\braceld$}%
  \braceld\leaders\vrule \@height\ht\z@ \@depth\z@\hfill\braceru$}

\makeatother

\usepackage{amscd, amsfonts, amsmath,amssymb,amsthm, mathtools} 
\usepackage{tabularx,ltablex}
\usepackage{stmaryrd}
\usepackage[all,color]{xy}
\usepackage[hidelinks]{hyperref}
\usepackage{cleveref}
\usepackage{float}
\usepackage{changes}
\usepackage{mathrsfs} 
\usepackage{amsmath,latexsym, bbold}
\usepackage{endnotes}
\usepackage{ulem}
\usepackage{aligned-overset}
\usepackage[margin=1in]{geometry}
\usepackage{enumitem} 
\allowdisplaybreaks
\newtheorem{introques}{Question}

\newtheorem{introdethm}[introques]{Theorem}

\newtheorem{thm}{Theorem}[subsection]
\newtheorem*{maintheorem*}{Main Theorem}
\newtheorem{lemma}[thm]{Lemma} 
\newtheorem{deflem}[thm]{Definition-Lemma} 
\newtheorem{proposition}[thm]{Proposition} 
\newtheorem{Cor}[thm]{Corollary}

\theoremstyle{definition}
\newtheorem{defn}[thm]{Definition}
\newtheorem{Ex}[thm]{Example} 
\newtheorem{notation}[thm]{Notation} 
\newtheorem{remark}[thm]{Remark} 

\newcommand{\id}{\mathrm{id}}

\DeclareMathOperator{\PGL}{PGL}
\DeclareMathOperator{\GL}{GL}

\DeclareMathOperator{\Ext}{Ext}
\DeclareMathOperator{\Hom}{Hom}

\DeclareMathOperator{\comod}{comod}

\DeclareMathOperator{\Imm}{Im}

\DeclareMathOperator{\op}{op}

\DeclareMathOperator{\alt}{Alt}
\DeclareMathOperator{\Sym}{Sym}
\DeclareMathOperator{\fd}{fd}
\newcommand{\HKQS}{(H,K)\text{-QS}}
\newcommand{\GLQS}{\mathcal{GL}\text{-QS}}
\newcommand{\SLQS}{\mathcal{SL}\text{-QS}}

\newcommand{\uaut}{\underline{\rm aut}}

\newcommand\res[2]{{
  \left.\kern-\nulldelimiterspace 
  #1 
  \littletaller 
  \right|_{#2} 
  }}

\newcommand{\littletaller}{\mathchoice{\vphantom{\big|}}{}{}{}}

\usepackage{tikz}
\usetikzlibrary{matrix}
\usepackage{graphicx,ctable,booktabs}
\usetikzlibrary{arrows.meta}
\usepackage{tikz-cd}

\usepackage{multirow}

\def\ob{\operatorname {ob}}

\def\bdt{\Delta}

\def\ot{\otimes}

\newcommand{\kk}{\Bbbk}

\newcommand{\mc}{\mathcal}

\newcommand{\cotens}{\mathbin{\square}}

\def\ra{\rightarrow}

\def\t{\text}
\def\it{\textit}

\def\vps{\varepsilon}

\numberwithin{equation}{section}

\title{Quantum-symmetric equivalence for superpotential algebras}

\author[Huang]{Hongdi Huang}
\address{(Huang) Department of Mathematics, Shanghai University, Shanghai 200444, China.}
\address{(Huang) Newtouch Center for Mathematics of Shanghai University, Shanghai 200444, China.}

\email{hdhuang@shu.edu.cn}
\author[Nguyen]{Van C. Nguyen}
\address{(Nguyen) Department of Mathematics, United States Naval Academy, Annapolis, MD 21402, U.S.A.}
\email{vnguyen@usna.edu}

\author[Vashaw]{Kent B. Vashaw}
\address{(Vashaw) Department of Mathematics and Statistics,
University at Albany SUNY,
Albany, NY 12222, U.S.A.}
\email{kvashaw@albany.edu}

\author[Veerapen]{Padmini Veerapen}
\address{(Veerapen) Department of Mathematics, Indiana University Northwest, Gary, IN 46408, U.S.A.}
\email{pveerap@iu.edu}

\author[Wang]{Xingting Wang}
\address{(Wang) Department of Mathematics, Louisiana State University, Baton Rouge, Louisiana 70803, U.S.A.}
\email{xingtingwang@math.lsu.edu}

\date\today

\keywords{quantum-symmetric equivalence, superpotential algebra, Artin--Schelter regular, bi-Galois object, cogroupoid}

\begin{document}
\maketitle

\begin{abstract}
For two connected graded algebras, we introduce the notion of quantum-symmetric equivalence defined relative to two fixed Hopf coactions. In particular, we consider algebras associated to nondegenerate twisted superpotentials and investigate the $\mathcal{GL}$-type and $\mathcal{SL}$-type quantum-symmetric equivalences using Bichon's reformulation of bi-Galois objects in the language of cogroupoids. For the $\mathcal{GL}$-type, we characterize the Koszul Artin--Schelter regularity, or equivalently, the Koszul twisted Calabi--Yau property, of a superpotential algebra as the non-vanishing of the bi-Galois object in the associated cogroupoid. For the $\mathcal{SL}$-type, we apply the pivotal structure of the comodule categories to study numerical invariants for $\mathcal{SL}$ quantum-symmetric equivalence, including the quantum Hilbert series of the superpotential algebra.
\end{abstract}

\section{Introduction}
The notion of superpotentials has proven to be a useful way of describing coordinate rings of noncommutative projective spaces. For example, Bocklandt \cite{BCY3} found that any graded Calabi--Yau (CY) algebra of dimension $3$ can be represented as a path algebra with relations derived from a superpotential. Dubois-Violette \cite{Dubois-Violette2007} and Bocklandt, Schedler and Wemyss \cite{BSW} generalized this result to any $N$-Koszul, (twisted) CY algebra. A particularly interesting question is to determine which superpotentials give rise to superpotential algebras which are (twisted) CY algebras. For connected graded algebras, the work by Reyes, Rogalski, and Zhang in \cite{RRZ2014} leads to a similar question about which superpotential algebras are Artin--Schelter (AS) regular algebras. By the work of  by Artin, Tate, and Van den Bergh in the late 1980s, AS-regularity is an essential concept in the development of noncommutative projective algebraic geometry. In their influential papers \cite{AS1987, ATV1990, ATV1991}, Artin, Schelter, Tate, and Van den Bergh classified AS-regular algebras in dimension $3$. The classic algebro-geometric methods developed in lower dimensions seem not to work effectively in dimensions $4$ and higher. Hence, the classification of higher-dimensional AS-regular algebras remains elusive. Recent work of Bhatoy, Ingalls, LaRoche, and Nookala suggests that viewing AS-regular algebras from the superpotential perspective could be fruitful \cite{BI2026,BILN2026}. 

In this work, we investigate superpotential algebras through their quantum symmetries, which are encoded in various universal quantum groups coacting on them. Quantum symmetries are crucial in the realm of noncommutative projective algebraic geometry (see \cite{Manin2018}), and Hopf (co)actions on noncommutative graded algebras, particularly on AS-regular algebras, is an extensive active research field (see e.g., \cite{AGV2023, CKWZ2018, CKWZ2019, CS2019, RVdB2017}). In previous joint work with Ure \cite{HNUVVW2024-2}, the authors introduced the concept of quantum-symmetric equivalence, describing when two connected graded algebras share the same quantum symmetries (which relates to the coactions of their associated universal quantum groups, in the sense of Manin), and proved that Koszul AS-regular algebras with the same global dimension form a single quantum-symmetric equivalence class. For superpotential algebras, two notable types of universal quantum group coactions exist: the \(\mathcal{GL}\)-type, which resembles Manin's universal quantum group \cite{Manin2018}, and the \(\mathcal{SL}\)-type, which can be regarded as universal cosovereign quantum groups \cite{Bichon2001}. We provide a criterion for when two superpotential algebras are $\mathcal{GL}$-type and $\mathcal{SL}$-type quantum-symmetrically equivalent in this setting. Our work builds on Bichon's reformulation of bi-Galois objects using the language of cogroupoids \cite{Bichon2014} by explicitly applying the \(\mathcal{GL}\) and \(\mathcal{SL}\) cogroupoids, which are constructed from nondegenerate twisted superpotentials \cite{HNUVVW2024-3}. An application of our results is the use of the \(\mathcal{GL}\)-type cogroupoid  to connect the classification of AS-regular algebras to the non-vanishing of specific bi-Galois objects. 
Meanwhile, the \(\mathcal{SL}\)-type provides numerical invariants derived from the pivotal structure on the comodule categories of the \(\mathcal{SL}\)-type universal quantum groups, which helps us to distinguish superpotential algebras that exhibit different quantum-symmetric behaviors. Our approach offers an alternative method for classifying AS-regular algebras by examining the connected components of the cogroupoid formed by their associated superpotentials from a quantum-symmetric perspective. It is important to note that while Artin and Schelter included finite Gelfand--Kirillov (GK) dimension in their definition of AS-regularity, we do not impose this requirement in this paper; our approach can be applied to a wide range of Koszul AS-regular rings over an arbitrary base field that may have infinite GK dimensions.

The paper is organized as follows. In \Cref{sec:prelim}, we recall necessary background on AS-regularity, bi-Galois objects, and the cogroupoids $\mathcal{GL}_m$ and $\mathcal{SL}_m$ in the context of superpotential algebras. \Cref{sec:QSsuperpotential} contains the main results of the paper. In particular, we extend the notion of quantum-symmetric (QS) equivalence between two connected graded algebras $A,B$ using their associated Manin's universal quantum groups (cf.~\cite[Definition A]{HNUVVW2024-2}) to an equivalence, called $\HKQS$ equivalence, using two Hopf algebras $H,K$ respectively coacting on them (see \Cref{def:HKQS}). When $A,B$ are locally finite, we show that $\HKQS$ equivalence implies QS equivalence.

\begin{introdethm}[\Cref{RQS}]
\label{thmA}
    Let $A$ and $B$ be two connected graded algebras that are locally finite. If $A$ and $B$ are $\HKQS$ equivalent, for some Hopf algebras $H$ and $K$ right coacting on $A$ and $B$ respectively, preserving their grading, then $A$ and $B$ are QS equivalent. 
\end{introdethm}

When $A=A(e,N)$ and $B=A(f,N)$ are superpotential algebras associated to two nondegenerate twisted superpotentials $e \in U^{\otimes m}$ and $f \in V^{\otimes m}$ respectively, for integers $2 \leq N \leq m$ (see \Cref{defn:superpotential}), we use the language of cogroupoids and specialize $\HKQS$ equivalence to $H=\mathcal{GL}_m(e)$ and $K=\mathcal{GL}_m(f)$ (or respectively to $H=\mathcal{SL}_m(e)$ and $K=\mathcal{SL}_m(f)$) to introduce the notion of $\GLQS$ equivalence (or respectively $\SLQS$ equivalence), cf.~\Cref{GLQSandSLQS}. We show that in this setting, the $\mathcal{GL}$-type (or respectively $\mathcal{SL}$-type) quantum-symmetric equivalence can be determined by the non-vanishing of the bi-Galois object $\mathcal{GL}_m(e,f)$ (or respectively $\mathcal{SL}_m(e,f)$). 

\begin{introdethm}[\Cref{SQSMTE} and \Cref{thm:SLQS}]
    Let $2 \leq N \leq m$ be integers, and $A(e,N)$ and $A(f,N)$ be two superpotential agebras associated to some nondegenerate twisted superpotentials $e\in U^{\otimes m}$ and $f\in V^{\otimes m}$, respectively. 
    \begin{itemize}
        \item[(1)] If $\mathcal{GL}_m(e,f)\neq 0$, then $A(e,N)$ and $A(f,N)$ are $\GLQS$ equivalent;
        \item[(2)] If $\mathcal{SL}_m(e,f)\neq 0$, then $A(e,N)$ and $A(f,N)$ are $\SLQS$ equivalent.
    \end{itemize} 
In either case, $A(e,N)$ and $A(f,N)$ are QS equivalent by \Cref{thmA}.  
\end{introdethm}
Under an additional assumption that $e,f$ are $L$-tracable for some $2 \leq L \leq N \leq m$ (cf.~\Cref{defn:traceable}), the converse statements of the above result also hold. 
This non-vanishing property of the bi-Galois object can further characterize the Koszul AS-regular property of the superpotential algebra as follows. 

\begin{introdethm}[\Cref{ASSuper} and \Cref{AS2cocyle}]
    Let $m\ge 2$ be any integer and $f\in V^{\otimes m}$ be a nondegenerate twisted superpotential. The following are equivalent.
    \begin{enumerate}
        \item[(1)] The associated superpotential algebra $A(f,2)$ is Koszul AS-regular; 
        \item[(2)] $\mathcal{GL}_m(e,f)\neq 0$ for all nondegenerate twisted superpotentials $e \in U^{\otimes m}$ such that $A(e,2)$ is Koszul AS-regular; 
        \item[(3)] $\mathcal{GL}_m(e,f)\neq 0$ for some nondegenerate twisted superpotential $e \in U^{\otimes m}$ such that $A(e,2)$ is Koszul AS-regular. 
    \end{enumerate}
Consequently, let $A$ be a Koszul AS-regular algebra of finite global dimension and $B$ any connected graded algebra finitely generated in degree one. Then $A$ and $B$ are QS equivalent if and only if $B$ is a Koszul AS-regular algebra of the same global dimension as $A$. In this case, $A$ and $B$ are 2-cocycle twists of each other if and only if they have the same Hilbert series.     
\end{introdethm}
The results in Theorem C shows that all superpotentials giving Koszul AS-regular algebras of dimension $m$ form a connected component of the cogroupoid $\mathcal{GL}_m$. We then characterize this connected component using several equivalent descriptions (see \Cref{IsoAS}). 

In \Cref{subsec:SLQS}, using the pivotal structure of $\mathcal{SL}_m$, we discuss some properties of the quantum dimension of any finite-dimensional $\mathcal{SL}_m(e)$-comodule and the quantum Hilbert series $q_{A(e,N)}(t)$ of the superpotential algebra $A(e,N)$ (cf.~\Cref{defn:qdim}). In particular, we show in \Cref{prop:SameH} that if $\mathcal{SL}_m(e,f) \neq 0$ then the two superpotential algebras $A(e,N)$ and $A(f,N)$ share the same quantum Hilbert series $q_{A(e,N)}(t)=q_{A(f,N)}(t)$, proving once again the usefulness of the non-vanishing of this bi-Galois object. When $N=m=2$, we give several characterizations for the $\SLQS$ equivalence, in relation to the non-vanishing of $\mathcal{SL}_m(e,f)$, quantum dimension, and quantum Hilbert series (see \Cref{equivm2}). 

In \Cref{sec:Example}, we illustrate our results by using the computer algebra system Magma to test the non-vanishing of $\mathcal{GL}_{3}(e,f)$ for some examples of non-noetherian algebras of dimension $3$, with $4$ generators that are given by nondegenerate superpotentials coming from the symmetrization of cubic surfaces in $\mathbb P^3$. We prove that several superpotential algebras arising in this way are not AS-regular, and we conversely conjecture, based on the non-termination of Magma in computing a Groebner basis for $\mathcal{GL}_3(e,f)$, that several large families are AS-regular (see \Cref{tab:Singular cubic}).

Finally, in the following diagram, for any connected graded algebras $A,B$ and some Hopf algebras $H,K$ respectively coacting on them preserving their grading, we visually summarize the connection between various quantum-symmetric equivalences between $A$ and $B$ that are studied in this paper. The dashed arrows are non-trivial implications due to the results discussed in this paper:

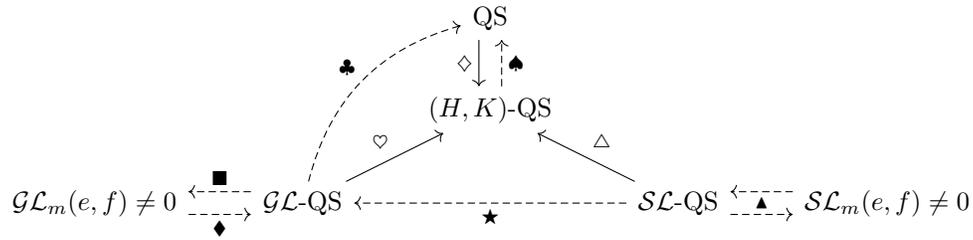
\begin{center}
\begin{tikzcd}
&& \textnormal{QS} \arrow[d, shift left=-1ex, "\diamondsuit" swap] & &\\
&& \HKQS \arrow[u, dashed, shift right=1ex, "\spadesuit" swap] & & \\
{\mathcal{GL}_m(e,f)\neq 0} \arrow[r, dashed, shift left=-1ex, "\blacklozenge" swap]  & \GLQS \arrow[l, dashed, shift right=1ex, "\blacksquare" swap] \arrow[ruu, "\clubsuit", dashed, bend left] \arrow[ru, "\heartsuit"] && \SLQS  \arrow[r, dashed, shift left=-1ex, "\blacktriangle"] \arrow[ll, "\bigstar", dashed] \arrow[ul, "\triangle" swap]& {\mathcal{SL}_m(e,f)\neq 0} \arrow[l, dashed, shift left=-1ex]
\end{tikzcd}
\end{center}
where 
\begin{itemize}[label={}]
\item $\diamondsuit$: when $H=\underline{\rm aut}^r(A)$ and $K=\underline{\rm aut}^r(B)$ are Manin's universal quantum groups (cf.~\cite[Definition 2.1.1]{HNUVVW2024-2}) associated to $A,B$ respectively;
\smallskip

\item $\heartsuit$: when $A=A(e, N)$ and $B=A(f, N)$ are superpotential algebras corresponding to nondegenerate twisted superpotentials $e,f$ respectively (cf.~\Cref{defn:superpotential}); $H=\mathcal{GL}_m(e)$ and $K=\mathcal{GL}_m(f)$ are Hopf algebras defined in \Cref{cogp} for any integers $2 \leq N \leq m$;
\smallskip

\item $\triangle$: when $A=A(e, N)$ and $B=A(f, N)$; $H=\mathcal{SL}_m(e)$ and $K=\mathcal{SL}_m(f)$ are defined in \Cref{cogp};
\smallskip

\item $\spadesuit$: by \Cref{RQS}, when $A, B$ are locally finite;
\smallskip

\item $\clubsuit$: by \Cref{SQSAS}, when $A=A(e,N)$ and $B=A(f,N)$ for $L$-traceable $e,f$ (cf.~\Cref{defn:traceable}) for some $2 \leq L \leq N \leq m$; $H=\mathcal{GL}_m(e)\cong \underline{\rm aut}^r(A(e, N))$ and $K=\mathcal{GL}_m(f)\cong \underline{\rm aut}^r(A(f,N))$;
\smallskip

\item $\blacklozenge$: by \Cref{SQSMTE}, when $A=A(e, N)$ and $B=A(f, N)$ and $\mathcal{GL}_m(e,f)$ is defined in \eqref{eq:alg};
\smallskip

\item $\blacksquare$: by \Cref{SQSBG}, when $A=A(e, N)$ and $B=A(f, N)$ for $L$-traceable $e,f$ with $2 \leq L \leq N \leq m$;\smallskip

\item $\bigstar$: by \Cref{thm:SLQS} and \Cref{SQSMTE}, due to the fact that $\mathcal{SL}_m(e)$ and $\mathcal{SL}_m(f)$ are respective quotients of $\mathcal{GL}_m(e)$ and $\mathcal{GL}_m(f)$;
\smallskip

\item $\blacktriangle$: by \Cref{thm:SLQS}, when $A=A(e, N)$ and $B=A(f, N)$. For the direction $\SLQS$ implies $\mathcal{SL}_m(e,f)\neq 0$, we need to assume $e,f$ are $L$-traceable for some $2 \leq L \leq N \leq m$.
\end{itemize}

\subsection*{Conventions} 

Throughout, let $\kk$ be a base field with $\otimes$ taken over $\kk$ unless stated otherwise. For any $\kk$-vector space $V$, we denote $V^* \coloneqq \Hom_\kk(V,\kk)$ its $\kk$-dual. A $\mathbb Z$-graded algebra $A=\bigoplus_{i\in \mathbb Z} A_i$ is called \emph{connected graded} if $A_i=0$ for $i<0$ and $A_0=\kk$. We use the Sweedler notation for the coproduct in a coalgebra $B$: for any $h \in B$, $\Delta(h) = \sum h_1 \otimes h_2 \in B \otimes B$. The category of right $B$-comodules is denoted by ${\rm comod}(B)$, and  ${\rm comod}_{\fd}(B)$ denotes its full subcategory consisting of finite-dimensional right $B$-comodules.

\subsection*{Funding} 
This work was supported by the National Science Foundation [DMS-2201146 to V.C.N., DMS-2103272 to K.B.V.]; and an AMS-Simons Travel Grant [to K.B.V.].

\subsection*{Acknowledgement} 
 The authors thank the referee for their helpful suggestions which improved the paper. This material is based upon work supported by the SQuaRE program at the American Institute of Mathematics (AIM) in Pasadena, California in June 2024; by the National Science Foundation Grant No. DMS-1928930 while the authors were in residence at the Simons Laufer Mathematical Sciences Institute (SLMath) in Berkeley, California in July 2024; and by the National Science Foundation Grant No. DMS-2201146 while Vashaw, Veerapen and Wang visited Nguyen at the U.S. Naval Academy in June 2025. The authors thank AIM, SLMath, and the U.S. Naval Academy for providing a supportive and mathematically rich environment. 

\medskip 

\noindent
{\bf Disclaimer.} The views expressed in this article are those of the authors and do not reflect the official policy or position of the U.S. Naval Academy, Department of the Navy, the Department of Defense, or the U.S. Government.


\section{Preliminaries}
\label{sec:prelim}

\subsection{Artin--Schelter regular algebras and Calabi--Yau algebras}

AS-regular algebras, originally introduced in \cite{AS1987}, are viewed as noncommutative analogs of polynomial rings. Below, we review the classical definition of AS-regularity. 

\begin{defn}[{\cite[Definition 2.1]{Rogalski2016}}]
\label{defn:AS-regular}
A connected graded algebra $A$ is called \emph{Artin--Schelter regular} (AS-regular) of dimension $d$ if the following conditions hold.
\begin{itemize}
    \item[(AS1)] $A$ has finite global dimension $d$, and
    \item[(AS2)] $A$ is \emph{Gorenstein}, that is, for some integer $l$,
\[
\Ext_A^i(\kk,A)=\begin{cases}
\kk(l) & \text{if}\ i=d\\
0 & \text{if}\ i\neq d,  
\end{cases}
\]
where $\kk$ is the trivial module $A/A_{\ge 1}$. We call the above integer $l$ the \emph{AS index of $A$}.
\end{itemize} 
\end{defn}

\begin{remark}\label{infinite GK}
We note that the original definition in \cite{AS1987,ATV1990} also assumes that $A$ has polynomial growth or, equivalently, that $A$ has finite Gelfand--Kirillov (GK) dimension. However, polynomial growth is not invariant under the Morita--Takeuchi equivalence. Hence, we follow the conventions of e.g. \cite{MS2016, RVdB2017, Zhang1998} to explore a larger family of noncommutative algebras, omitting the finite GK dimension assumption in this paper. 
\end{remark}

In a broader sense, AS-regular algebras represent a distinct family of connected graded twisted Calabi--Yau algebras, which generalize Ginzburg’s Calabi--Yau algebras. Given an algebra $A$, we write $A^{\op}$ for the algebra $A$ with opposite multiplication and $A^e \coloneqq A\otimes A^{\op}$ for the enveloping algebra of $A$. Hence, an $A$-bimodule can be identified with a left $A^e$-module, that is, an object in ${\rm mod}(A^e)$. Let $M$ be an $A$-bimodule, and $\mu,\nu$ be algebra automorphisms of $A$. Then $\!^\nu M^\mu$ denotes the induced $A$-bimodule such that $\!^\nu M^\mu=M$ as a $\kk$-vector space, with the left and right actions of $A$ given by $\mu,\nu$ via
\[a\,\cdot_\nu\, m\, \cdot_\mu\, b = \nu(a)m\mu(b)\]
for all $a,b\in A$ and $m\in M$. In particular, we will omit $\mu$ or $\nu$ if it is the identity.

\begin{defn}[{\cite[Definition 9.3]{Rogalski2016}}]\label{defn tcy}
An algebra $A$ is called  a \emph{twisted Calabi--Yau (CY) algebra  of dimension
$d$} if
\begin{enumerate}
\item[(CY1)] $A$ is \emph{homologically smooth}, that is, $A$ has
a length $d$ resolution by finitely generated projective
$A^e$-modules; and 
\item[(CY2)] There is an automorphism $\mu$ of $A$ such that
\begin{equation*}
\Ext_{A^e}^i(A,A^e)\cong
\begin{cases}
A^\mu,&i=d\\
0,& i\neq d
\end{cases}
\end{equation*}
as $A^e$-modules.
\end{enumerate}
\end{defn}

In the above definition of a twisted CY algebra $A$, we call the associated automorphism $\mu$ the \it{Nakayama automorphism} of $A$, which is unique up to conjugation by an inner automorphism. Consequently, a CY algebra in the sense of Ginzburg is a twisted CY algebra whose Nakayama automorphism is an inner automorphism. In particular, when $A$ is connected graded, we can require its Nakayama automorphism to be a graded algebra automorphism, in which case it is uniquely determined. 

The following proposition indicates that the twisted CY and AS-regular properties coincide for connected graded algebras. Thus, as we focus on connected graded algebras in this paper, any result for AS-regularity also implies for twisted CY property.

\begin{proposition}[{\cite[Lemma 1.2]{RRZ2014}}]
    Let $A$ be a connected graded algebra. Then $A$ is twisted CY if and only if it is AS-regular.
\end{proposition}

\subsection{Nondegenerate twisted superpotentials} 
\label{subsect:nondegentwistsuppot}

\begin{defn}
\label{defn:superpotential}
Let $2 \leq N \leq m$ be integers and $V$ an $q$-dimensional $\kk$-vector space. Let $\phi: V^{\otimes m} \rightarrow V^{\otimes m}$ be the linear map defined by 
\[\phi\left( v_1 \otimes \cdots \otimes v_m \right) \coloneqq  v_m \otimes v_1 \otimes \cdots \otimes v_{m-1}, \quad \text{for any $v_i \in V$}.\]
\begin{enumerate}
\item[(1)] Take an element $s \in V^{\otimes m}$. 
\begin{itemize}
\item[(i)] We say $s$ is a \emph{twisted superpotential} if there is $\mathbb{P} \in \GL(V)$ so that \[\left( \mathbb{P} \otimes \id^{\otimes(m-1)} \right) \phi(s) = s.\]
\item[(ii)] A twisted superpotential $s$ is said to be 
\emph{nondegenerate} if for $\nu_1\in V^*$, $(\nu _1\otimes \nu _2\otimes \cdots \otimes \nu _m)(s)=0$ for all $\nu_2,\ldots,\nu_m\in V^*$ implies that $\nu _1=0$. 
\end{itemize}
\item[(2)] Given a twisted superpotential $s \in V^{\otimes m}$, the \emph{superpotential algebra} associated to $s$ is defined as
\[A(s,N) \coloneqq TV/\langle\partial^{m-N} (\Bbbk s) \rangle,\]
where $TV$ is the tensor algebra on $V$,
\begin{align*} 
\partial  (\Bbbk s) &= \left\{ \left( \nu \otimes \id^{\otimes (m-1)}\right) (\alpha s) \mid \nu \in V^*, \alpha \in \Bbbk \right\}, 
\text{ and } \\
\partial^{i+1} (\Bbbk s) &= \partial \left( \partial^i ( \Bbbk s)\right), \quad \text{for all }i \geq 0. \end{align*} 
    \end{enumerate}
\end{defn} 

In the above definition, $s\in V^{\otimes m}$ is called a \emph{superpotential} if $\phi(s)=s$, that is, $\mathbb P=\mathbb I$ the identity matrix. Without loss of generality, we may consider a twisted superpotential as a linear map $s: \kk \to V^{\otimes m}$. We use this identification interchangeably throughout the paper. Moreover, nondegenerate twisted superpotentials are closely related to preregular forms (see e.g., \cite{Dubois-Violette2007}). 

Recall that if $N \ge 2$ is an integer, a connected graded algebra $A$ finitely generated in degree 1 with $N$-homogeneous relations is called $N$-Koszul if  $\Ext^i_A(\kk,\kk)$ is concentrated in a single degree for all $i$. $N$-Koszul algebras were introduced by Berger in \cite{Berger2001}, and generalize the usual notion of Koszul algebras; in particular, when we write Koszul, this is the same as 2-Koszul under Berger's definition. 

The combined result due to Dubois-Violette and Mori--Smith demonstrates that every $N$-Koszul AS-regular algebra arises from superpotential algebras, and it also identifies which superpotential algebras yield CY algebras. 
 
\begin{thm}[{\cite[Theorem 11]{Dubois-Violette2007}, \cite[Corollary 4.5]{MS2016}}]
\label{NKAS}
Let $A$ be an $N$-Koszul, Artin--Schelter regular algebra of global dimension $d$. There is a unique (up to multiplication by a nonzero scalar) nondegenerate twisted superpotential $s\in A_1^{\otimes m}$ with integers $2 \leq N \leq m$ such that $A=A(s,N)$. Moreover:
\begin{itemize}
    \item[(1)] $m$ is the AS index for $A$. When $N=2$, $m=d$. When $N\ge 3$, $d$ is odd and $m=\frac{1}{2}N(d-1)+1$.
    \item[(2)] The Nakayama automorphism $\mu$ of $A(s,N)$ is determined by $\mu|_{A_1}=(-1)^{d+1}(\mathbb P^T)^{-1}$, where $(s,\mathbb P)$ is the corresponding nondegenerate twisted superpotential in $A_1^{\otimes m}$.
    \item[(3)] $A(s,N)$ is CY if and only if $\mathbb P=(-1)^{d+1}\mathbb I$ if and only if $\phi(s)=(-1)^{d+1}s$ where $\phi: A_1^{\otimes m}\to A_1^{\otimes m}$ is the cyclic permutation. 
\end{itemize}
\end{thm}

Let $2\leq N\leq m$ be integers and $V$ any $q$-dimensional vector space. We identify $V^{\otimes m}$ as an $mq$-dimensional vector space and define ${\rm Gr}(1,mq)$ to be the Grassmannian parameterizing the set of all one-dimensional subspaces of $V^{\otimes m}$. Let $f\in V^{\otimes m}$ be a nondegenerate twisted superpotential. Then $\kk f$ can be identified as an element in ${\rm Gr}(1,mq)$. Now, we define a subset of the Grassmannian ${\rm Gr}(1,mq)$ as 
\begin{equation}
\label{ASclass}
   {\rm AS}(m,N,q) \coloneqq \{\kk f\in {\rm Gr}(1,mq)\,|\, \text{$A(f,N)$ is $N$-Koszul AS-regular}\}\subseteq {\rm Gr}(1,mq). 
\end{equation}
It is clear that ${\rm Gr}(1,mq)$ has a natural $\PGL_q(\kk)$ action by acting on $V$ diagonally in $V^{\otimes m}$, which makes ${\rm AS}(m,N,q)$ a $\PGL_q(\kk)$-invariant subset.  There is a subset of ${\rm AS}(m,N,q)$ corresponding to those $N$-Koszul AS-regular algebras that are Calabi--Yau algebras of dimension $d=2(m-1)/N+1$, which we denote by 
\begin{equation}
\label{CYclass}
    {\rm CY}(m,N,q) \coloneqq \{\kk f\in {\rm AS}(m,N,q)\,|\, \phi(f)=(-1)^{d+1}f\}\subseteq {\rm AS}(m,N,q),
\end{equation}
where $\phi: V^{\otimes m}\to V^{\otimes m}$ is the cyclic permutation. For short, we refer to any $\kk f \in {\rm AS}(m,N,q)$ or ${\rm CY}(m,N,q)$ as an $N$-Koszul AS-regular superpotential or an $N$-Koszul Calabi--Yau superpotential, respectively.

\subsection{Bi-Galois objects}

We provide some background on Galois objects, also referred to as Galois extensions in the literature, resulting from Hopf coactions. Such objects stem from classical Galois correspondences for field extensions (see e.g., \cite{Montgomery2009} for a survey). In this work, we examine Schauenburg's \cite{Sch1996} generalization of the construction of $H$-Galois objects for a Hopf algebra $H$, as studied in, e.g., \cite{Ulbrich1987,VOZ1994}. 

\begin{defn}[{\cite{Sch1996}}] 
\label{defn:bi-Galois}
Let $H$ and $K$ be two Hopf algebras. A \emph{left $H$-Galois object} is a left $H$-comodule algebra $T$ such that the linear map
\[
T\otimes T\xrightarrow{\rho\otimes \id} H\otimes T\otimes T\xrightarrow{\id \otimes m} H\otimes T
\]
is bijective, where $\rho: T\to H\otimes T$ and  $m: T\otimes T\to T$ are the corresponding left $H$-comodule algebra structure maps on $T$. Two left $H$-Galois objects are \emph{isomorphic} if they are isomorphic as $H$-comodule algebras. A \emph{right $K$-Galois object} can be defined analogously, and an \emph{$H$-$K$-bi-Galois object} is an $H$-$K$-bicomodule algebra that is both a left $H$-Galois object and a right $K$-Galois object. An $H$-$K$-bi-Galois object $T$ is called \emph{cleft} if $T$ is isomorphic to $H$ as a left $H$-comodule, and isomorphic to $K$ as a right $K$-comodule. 
\end{defn}

We say that two Hopf algebras $H$ and $K$ are \emph{Morita--Takeuchi equivalent} if there is some monoidal equivalence described by $(F,\xi)$, where   \[F: {\rm comod}(H)\to{\rm comod}(K)
\]
is an equivalence of categories, together with a natural isomorphism 
\begin{align*}
 \xi_{W,M}: F(W)\otimes F(M)~\cong~F(W\otimes M),
\end{align*}
for any two $H$-comodules $W$ and $M$ satisfying some compatible conditions in \cite[Definition 2.4.1]{EGNO2015}. For any tensor product of $H$-comodules $W_1\otimes \cdots \otimes W_n$, one can construct a natural isomorphism between $F(W_1)\otimes \cdots \otimes F(W_n)$ and $F(W_1\otimes \cdots \otimes W_n)$ by composing several instances
of the tensor products of $\xi$ and the identity. Any of the above isomorphisms is unique by the compatible condition on the monoidal structure and will be denoted by $\xi: F(W_1)\otimes \cdots \otimes F(W_n)\to F(W_1\otimes \cdots \otimes W_n)$ for simplicity. The \emph{duality transformation} of $(F,\xi)$ is the canonical isomorphism between $F(W^*)$ and $F(W)^*$, denoted by $\widetilde{\xi}: F(W^*)\to F(W)^*$, satisfying the following commutative diagrams: 
\begin{align}
\label{xitilde}
\xymatrix{
F(W^*)\otimes F(W)\ar[r]^{\xi}\ar[d]_-{\widetilde{\xi}\otimes \id} & F(W^*\otimes W)\ar[d]^-{F({\rm ev})}\\
F(W)^*\otimes F(W)\ar[r]^-{{\rm ev}} & \kk
}
\qquad \text{and} \qquad 
\xymatrix{
\kk\ar[r]^-{{\rm coev}}\ar[d]_-{F({\rm coev})} & F(W)\otimes F(W)^*\ar[d]^-{\id\otimes \widetilde{\xi}^{-1}}\\
F(W\otimes W^*)\ar[r]^-{\xi^{-1}} & F(W)\otimes F(W^*).
}
\end{align}
Moreover, $\widetilde{\xi}$ is uniquely determined by either of the commutative diagrams above (cf.~ \cite[\S 1]{Ng-Schauenburg}). By \cite[Proposition 1.16]{Bichon2014}, $(F,\xi)$ restricts to a monoidal equivalence between ${\rm comod}_{\rm fd}(H)$ and  ${\rm comod}_{\rm fd}(K)$. In particular, $F$ sends one-dimensional $H$-comodules to one-dimensional $K$-comodules. 

According to \cite{Sch1996}, any Morita--Takeuchi equivalence between $H$ and $K$ can be described in terms of bi-Galois objects (see \Cref{defn:bi-Galois}). That is, any monoidal equivalence $(F,\xi)$ can be explicitly given by the cotensor product with some $H$-$K$-bi-Galois object $T$:
\begin{equation}\label{eq:MTE} 
\begin{aligned}
F: {\rm comod}(H)~&\to~{\rm comod}(K) \\
W~&\mapsto~ W\cotens_HT, 
\end{aligned}
\end{equation}
where the cotensor product is given by 
\[
W\cotens_HT \coloneqq \{x\in W\otimes T\,|\, (\rho_W\otimes \id_T)(x)=(\id_W \otimes \rho_T)(x)\}.
\]
Here $\rho_W: W\to W\otimes H$ and $\rho_T: T\to H\otimes T$ are the corresponding $H$-comodule structure maps on $W$ and $T$, respectively. The right coaction of $K$ on $W\cotens_HT$ comes from that of $K$ on $T$, since $T$ is an $H$-$K$-bi-Galois object.
The natural isomorphisms giving $F$ the structure of a monoidal functor are given by (as in \cite[Proposition 1.4]{Ulbrich1987} for left comodules)
\begin{align*}
\xi_{W,M}: (W\cotens_HT)\otimes (M\cotens_HT)~&\cong~ (W\otimes M)\cotens_HT\\
\left(\sum w_i\otimes a_i\right)\otimes \left(\sum m_j\otimes b_j\right)~&\mapsto~ \sum w_i\otimes m_j\otimes a_ib_j,
\end{align*}
for any $w_i \in W, m_j \in M$, and $a_i,b_j \in T$. In the case that $T$ is a cleft bi-Galois object, the cotensor product $(-)\cotens_HT$ is given by a 2-cocycle \cite{DT86, BM1989}. For general background on 2-cocycles and 2-cocycle twists, see e.g., \cite{M2005}. 

The following result by Schauenburg \cite{Sch2004} shows that any left or right Galois object can be turned into a bi-Galois object satisfying a universal property. We will use this result several times later in the paper.

\begin{thm}[{\cite[Propositions 2.8.4 and 3.1.2, Corollary 3.1.3]{Sch2004}}]
\label{Schauenburg}
    Let $H$ be a Hopf algebra and $T$ be a left $H$-Galois object. There exists a Hopf algebra $K$ that right coacts on the algebra $T$ such that $T$ becomes an $H$-$K$-bi-Galois object. Moreover, $K$ is subject to the following universal property: for any Hopf algebra $L$ that right coacts on $T$, there is a unique bialgebra map $f: K\to L$ such that the diagram
    \[
    \xymatrix{
T\ar[rr]\ar[dr]&& T\otimes K\ar[dl]^-{\id \otimes f}\\
&T\otimes L & 
    }
    \]
    commutes. In particular, $K$ is unique up to isomorphism. 
\end{thm}

\subsection{Cogroupoids}
We now discuss the technical machinery necessary for investigating quantum symmetries related to superpotential algebras. This involves analyzing the Morita--Takeuchi equivalence between universal quantum groups associated to superpotential algebras, utilizing the concept of cogroupoids as introduced by Bichon \cite{Bichon2014}. 

\begin{defn}
\label{defn:cocategory}
A \emph{$\kk$-cocategory} $\mc{C}$ consists of:
\begin{enumerate}
\item A set of objects $\ob(\mc{C})$;
\item For any $X,Y\in \ob(\mc{C})$, a $\kk$-algebra $\mc{C}(X,Y)$;
\item For any $X,Y,Z\in \ob(\mc{C})$, $\kk$-algebra homomorphisms
\[ \bdt^Z_{XY}:\mc{C}(X,Y)\ra \mc{C}(X,Z)\ot \mc{C}(Z,Y) \qquad \t{and} \qquad \vps_X:\mc{C}(X,X)\ra \kk \]
such that for any $X,Y,Z,T\in \ob(\mc{C})$, the following diagrams commute
\[
\xymatrix{
\mc{C}(X,Y)\ar^-{\bdt^Z_{X,Y}}[rr]\ar_{\bdt^T_{X,Y}}[d] && \mc{C}(X,Z)\ot \mc{C}(Z,Y)\ar^-{\bdt^T_{X,Z}}[d]\\
\mc{C}(X,T)\ot \mc{C}(T,Y)\ar^-{\id\ot \bdt^Z_{T,Y}}[rr]&&\mc{C}(X,T)\ot \mc{C}(T,Z)\ot \mc{C}(Z,Y), } 
\]
\[\xymatrix
{\mc{C}(X,Y)\ar@{=}[rd]\ar_{\bdt^Y_{X,Y}}[d]&\\
\mc{C}(X,Y)\ot\mc{C}(Y,Y)\ar^-{\id \ot \vps_Y}[r]&\mc{C}(X,Y),} \qquad \qquad 
\xymatrix
{\mc{C}(X,Y)\ar@{=}[rd]\ar_{\bdt^X_{X,Y}}[d]&\\
\mc{C}(X,X)\ot\mc{C}(X,Y)\ar^-{\vps_X \ot \id }[r]&\mc{C}(X,Y).}
\]
\end{enumerate}
\end{defn}

For $a^{X,Y}\in \mc{C}(X,Y)$, we use Sweedler's notation to write 
\[ \bdt^Z_{X,Y}(a^{X,Y})=\sum a^{X,Z}_1\ot a^{Z,Y}_2. \]
From its definition, a cocategory with one object is just a bialgebra. In particular, $\mc{C}(X, X)$ is a bialgebra for any $X \in \ob(\mc{C})$. A cocategory $\mc{C}$ is said to be \emph{connected} if $\mc{C}(X,Y)$ is a nonzero algebra for any $X,Y\in \ob(\mc{C})$ \cite[Definition 2.4]{Bichon2014}.

\begin{defn}
\label{defn:cogroupoid}
A \emph{$\kk$-cogroupoid} $\mc{C}$ consists of a $\kk$-cocategory $\mc{C}$ together with linear maps
\[ S_{X,Y}:\mc{C}(X,Y)\longrightarrow \mc{C}(Y,X), \]
for any $X,Y\in \ob(\mc{C})$, such that the following diagram commutes
{\small \[\xymatrix{\mc{C}(Y,X) & \kk \ar[l]_-u & \mc{C}(X,X)\ar[d]_{\bdt_{X,X}^Y}\ar[r]^-{\vps_X}\ar[l]_-{\varepsilon_X} &\kk\ar[r]^-u&\mc{C}(X,Y)\\
\mc{C}(Y,X)\ot\mc{C}(Y,X) \ar[u]^m &&\mc{C}(X,Y)\ot\mc{C}(Y,X)\ar[rr]^{\id\ot S_{Y,X}} \ar[ll]_-{S_{X,Y} \otimes \id} &&\mc{C}(X,Y) \ot\mc{C}(X,Y). \ar[u]^m} \]} 
\end{defn}

In a cogroupoid $\mc{C}$, the bialgebra $\mc{C}(X,X)$ is a Hopf algebra for any $X\in \ob(\mc{C})$, with the antipode map given by $S_{X,X}$, and if $\mc{C}(X,Y)$ is nonzero then it is a $\mc{C}(X,X)$-$\mc{C}(Y,Y)$-bi-Galois object \cite[Proposition 2.8]{Bichon2014}. For other properties of cogroupoids, we refer the reader to \cite{Bichon2014}.

The following is Bichon's reformulation of Schauenburg's result \cite{Sch1996} about bi-Galois objects for Morita--Takeuchi equivalences in terms of cogroupoids. 

\begin{thm}[{\cite[Theorem 2.10]{Bichon2014}}]
Let $H$ and $K$ be Hopf algebras. Then the following conditions are equivalent.
 \begin{enumerate}
     \item[(1)] There exists a $\kk$-linear equivalence of monoidal categories ${\rm comod}(H)\to {\rm comod}(K)$;
     \item[(2)] There exists a connected cogroupoid $\mathcal C$ and two objects $X,Y \in \ob(\mc{C})$ such that $H \cong\mathcal C(X, X)$ and $K \cong\mathcal C(Y, Y )$.
 \end{enumerate}
\end{thm}

Moreover, according to \cite[Theorem 2.12]{Bichon2001}, in a cogroupoid $\mathcal C$, if $\mathcal C(X,Y)\neq 0$, then $\mathcal C(Y,X)\neq 0$ and we have $\kk$-linear equivalences of monoidal categories that are inverse of each other: 
\begin{align*}
\comod(\mathcal C(X,X))&\to \comod(\mathcal C(Y,Y)) & \comod(\mathcal C(Y,Y))&\to \comod(\mathcal C(X,X))\\
W&\mapsto W\cotens_{\mathcal C(X,X)}\mathcal{C}(X,Y) & M&\mapsto M\cotens_{\mathcal C(Y,Y)} \mathcal{C}(Y,X).
\end{align*}

\begin{defn}[{\cite[Definition 4.2.1]{HNUVVW2024-3}}]
A cogroupoid $\mathcal C$ is called \emph{copivotal}  
if for any $X\in \mathrm{ ob}(\mathcal C)$, there exists some  
character $\Phi_X: \mathcal C(X,X)\to \kk$ such that for any $X,Y\in \mathrm{ ob}(\mathcal C)$, 
\[
S_{Y,X}\circ S_{X,Y}=\Phi_X^{-1}*\id*\Phi_Y,
\]
where the right hand is defined as 
\begin{align*}
\Phi_X^{-1}*\id*\Phi_Y &:=\left(\mathcal C(X,Y)\xrightarrow{\Delta_{X,Y}^X} \mathcal C(X,X)\otimes \mathcal C(X,Y)\right.\\
&\quad \left.\xrightarrow{\id \otimes \Delta_{X,Y}^Y} \mathcal C(X,X)\otimes\mathcal C(X,Y)\otimes  C(Y,Y)\xrightarrow{\Phi_X^{-1}\otimes \id \otimes \Phi_Y} \mathcal C(X,Y)\right),
\end{align*}
and $\Phi_X^{-1} = \Phi_X \circ S_{X,X}$ is the convolution inverse of $\Phi_X$.
\end{defn}

\subsection{The cogroupoids $\mathcal{GL}_m$ and $\mathcal{SL}_m$}
\label{cogp}

We recall the cogroupoid $\mathcal{GL}_m$ constructed in \cite{HNUVVW2024-3} for the purpose of this paper and note the use of $m$-linear preregular forms in the original construction. We modify the initial construction to consider all nondegenerate twisted superpotentials. 

For any integer $m\geq 2$, $\mathcal{GL}_m$ is a $\kk$-cogroupoid, whose objects are all nondegenerate twisted superpotentials in $m$-fold tensor products of $\kk$-vector spaces and the Hom-space between two nondegenerate twisted superpotentials  $e\in U^{\otimes m}$ and $f\in V^{\otimes m}$, denoted by $\mathcal{GL}_m(e,f)$, is the $\kk$-algebra with $2(pq+1)$ generators 
    \[\mathbb A=(a_{ij})_{\substack{1\leq i\leq p \\ 1\leq j\leq q}}, \qquad \mathbb B=(b_{ij})_{\substack{1\leq i\leq q \\ 1\leq j\leq p}}, \qquad D^{\pm 1},\]
subject to the relations
\begin{equation}
\label{eq:alg}
\left. 
\begin{aligned}
    \sum_{1\leq i_1,\ldots,i_m\leq q}{f}_{i_1 \cdots i_m}a_{j_1i_1}\cdots a_{j_mi_m} &= {e}_{j_1\cdots j_m}D^{-1}, &\textnormal{ for any } 1\leq j_1,\dots, j_m\leq p,\\
    \sum_{1\leq i_1,\ldots,i_m\leq p}{e}_{i_1\cdots i_m}b_{j_mi_m}\cdots b_{j_1i_1} &= {f}_{j_1\cdots j_m}D,  &\textnormal{ for any } 1\leq j_1,\dots, j_m\leq q,\\
    DD^{-1} &= D^{-1}D=1, \text{ and } \\
    \mathbb{B}\mathbb{A} &= \mathbb{I}_{q\times q}.
\end{aligned}\right\}
\end{equation} 
Here we assume $\dim U=p$ and $\dim V=q$ with fixed bases $\{u_1,\ldots,u_p\}$ and $\{v_1,\ldots, v_q\}$, respectively.  We write  
\[
e=\sum_{1\leq i_1,\ldots,i_m\leq p} {e}_{i_1 \cdots i_m}\, u_{i_1}\otimes \dots \otimes u_{i_m}\quad \text{and}\quad 
f=\sum_{1\leq j_1,\ldots,j_m\leq q} {f}_{j_1 \cdots j_m}\, v_{j_1}\otimes \dots \otimes v_{j_m}
\]
for coefficients ${e}_{i_1 \cdots i_m}, {f}_{j_1 \cdots j_m}\in \kk$.

\begin{deflem}
[{\cite[Lemmas 3.1.3, 3.1.5, Definition 3.1.6]{HNUVVW2024-3}}]
\label{deflem}
For any integer $m\geq 2$, $\mathcal{GL}_m$ is a $\kk$-linear cogroupoid. For any vector spaces $U,V,W$ with $\dim U=p,$ $\dim V=q,$ and $\dim W=r$, and any nondegenerate twisted superpotentials $e\in U^{\otimes m}$, $f\in V^{\otimes m}$ and $g\in W^{\otimes m}$, the cogroupoid structure is given by the algebra maps
\[ \Delta = \Delta_{e,g}^f: \mathcal{GL}_m(e,g)\to \mathcal{GL}_m(e,f)\otimes \mathcal{GL}_m(f,g), \]
such that 
\begin{gather*}
\Delta(a_{ij}^{e,g})=\sum_{k=1}^q a^{e,f}_{ik}\otimes a^{f,g}_{kj},\quad \Delta(b_{ji}^{e,g})=\sum_{k=1}^q b^{e,f}_{ki}\otimes b^{f,g}_{jk}, \qquad\text{for } 1\leq i\leq p, 1\leq j\leq r,\\
\Delta((D^{e,g})^{\pm 1})=(D^{e,f})^{\pm 1}\otimes (D^{f,g})^{\pm 1}, 
\end{gather*} 
and
\[ \varepsilon_{e}: \mathcal{GL}_m(e)\to \kk,\]
such that $\varepsilon_{e}(a_{ij}^{e,e})=\varepsilon_{e}(b_{ji}^{e,e})=\delta_{ij}$, for $1\leq i,j\leq p$, and $\varepsilon_e((D^{e,e})^{\pm 1})=1$. 

Moreover, the algebra map
\[
S_{e,f}: \mathcal{GL}_m(e,f)\to \mathcal{GL}_m(f,e)^{\operatorname{op}}
\]
is defined by the formulas 
\begin{align*}
S_{e,f}\left(\mathbb A^{e,f}\right)&=\mathbb B^{f,e},\\
S_{e,f}\left(\mathbb B^{e,f}\right)&=\left(D^{f,e}\right)^{-1}\,\mathbb Q^{-T}\,\mathbb A^{f,e}\,\mathbb P^T\, D^{f,e},\\ S_{e,f}\left(\left(D^{e,f}\right)^{\pm 1}\right)&=\left(D^{f,e}\right)^{\mp 1},
\end{align*} 
where $\mathbb P\in \GL(U)$ and $\mathbb Q\in \GL(V)$ are invertible matrices associated with $e,f$ as in \Cref{defn:superpotential}.
\end{deflem}

We denote the generators of  $\mathcal{GL}_m(e,f)$ by $\mathbb A^{e,f}=(a_{ij}^{e,f}), \mathbb B^{e,f}=(b_{ij}^{e,f}),$ and $\left( D^{e,f}\right)^{\pm 1}$ when multiple nondegenerate twisted superpotentials are involved, and omit the superscripts when the context is clear. In particular, if $U=V$ and $e=f$, then we write $\mathcal{GL}_m (f)=\mathcal{GL}_m(f,f)$, which is the algebra $\mathcal{H}(f)$ associated to a preregular
form $f^*$ defined in \cite[Definition 5.17]{CWW2019}. In this paper we use the right comodule structure.

We recall the cogroupoid $\mathcal{SL}_m$ constructed from $\mathcal{GL}_m$ by factoring out by the ideal of relations given by setting all quantum determinants $D^{*,*}$ equal to $1$ \cite{HNUVVW2024-3}. In particular, we write $\mathcal{SL}_m(e)=\mathcal{SL}_m
(e,e)\coloneqq \mathcal{GL}(e,e)/(D^{e,e}-1)$, which is isomorphic to the universal quantum group defined by Bichon and Dubois-Violette \cite{BDV13} for a preregular form.   Moreover, the cogroupoid $\mathcal{SL}_m$ is copivotal by the following result.  

\begin{lemma}[{\cite[Lemma 4.2.5]{HNUVVW2024-3}}]\label{L:copi}
The cogroupoid $\mathcal{SL}_m$ is copivotal, with character $\Phi_e: \mathcal{SL}_m(e)\to \kk$ given by $\Phi_e(\mathbb A^{e,e})=\mathbb P^T$ and $\Phi_e(\mathbb B^{e,e})=\mathbb P^{-T}$, that is, 
\[
S_{f,e}\circ S_{e,f}=\Phi_e^{-1}* \id *\Phi_f
\]
in $\mathcal{SL}_m(e,f)$, for any two nondegenerate twisted superpotentials $e\in U^{\otimes m}$ and $f\in V^{\otimes m}$ with associated invertible matrices $\mathbb{P}$ and $\mathbb{Q}$, respectively.  
\end{lemma}

\section{Quantum-symmetric equivalence for superpotential algebras}
\label{sec:QSsuperpotential}

In this section, we apply the concept of quantum-symmetric equivalence to superpotential algebras and consequently, connect the AS-regularity of superpotential algebras with the Morita--Takeuchi equivalence of their associated universal quantum groups. We reformulate the classification problem of Koszul AS-regular algebras of global dimension $m$ in terms of connected components in the cogroupoid $\mathcal{GL}_m$ by classifying fiber functors from the comodule categories over the universal quantum groups of the polynomial algebras (see \Cref{ASSuper} and \Cref{IsoAS}).

\subsection{A characterization of quantum-symmetric equivalence}

Quantum-symmetric equivalence of connected graded algebras was originally introduced by Ure and the authors in \cite[Definition A]{HNUVVW2024-2}, motivated by connections between Zhang twists of graded algebras and 2-cocycle twists of Hopf algebras \cite{HNUVVW2024,HNVVW2025}, as well as the Raedschelders--Van den Bergh Theorem \cite{RVdB2017}. We will extend the notion of quantum-symmetric (QS) equivalence of connected graded algebras for Manin's universal quantum groups to QS equivalence for any Hopf coactions on the algebras. Let us first recall the definition of Manin's universal quantum groups \cite[Lemma 6.6]{Manin2018}; these universal quantum groups exist for connected graded $\kk$-algebras which are locally finite (that is, each graded component is finite-dimensional over $\kk$) \cite[Theorem 3.16 and Example 4.8(1)-(2)]{AGV2023}.

\begin{defn}[{\cite[Lemma 6.6]{Manin2018}}]
\label{def:ManinU} 
Let $A$ be any $\mathbb Z$-graded locally finite $\kk$-algebra. The \emph{right universal quantum group  $\uaut^r(A)$ associated to $A$} is the Hopf algebra that right coacts on $A$ preserving the grading of $A$ via $\rho: A\to A \otimes \uaut^r(A)$ satisfying the following universal property: if $H$ is any Hopf algebra that right coacts on $A$ preserving the grading of $A$ via $\tau: A\to A \otimes H$, then there is a unique Hopf algebra map $f: \uaut^r(A)\to H$ such that the diagram 
\begin{align}
\label{def:aut}
\xymatrix{
A\ar[r]^-{\rho}\ar[dr]_-{\tau} & A \otimes \uaut^r(A) \ar[d]^-{\id \otimes f} \\
& A \otimes H
}
\end{align}
commutes. 
\end{defn}

\begin{defn}[{\cite[Definition A]{HNUVVW2024-2}}] \label{DefnA}
Let $A$ and $B$ be two connected graded algebras finitely generated in degree one. We say $A$ and $B$ are \emph{quantum-symmetrically (QS) equivalent} if there is a monoidal equivalence between the comodule categories of their associated universal quantum groups 
\[
\comod(\uaut^r(A))~\overset{\otimes}{\cong}~\comod(\uaut^r(B))
\]
in the sense of Manin, where this equivalence sends $A$ to $B$ as comodule algebras. We denote the quantum-symmetric equivalence class of $A$ by $QS(A)$, which consists of all connected graded algebras that are quantum-symmetrically equivalent to $A$.
\end{defn}
We now generalize the notion of QS equivalence via replacing the coactions of Manin's universal quantum groups by any Hopf coactions.

\begin{defn}
\label{def:HKQS}
   Let $A$ and $B$ be any two connected graded algebras and let $H,K$ be two Hopf algebras coacting on $A,B$, respectively, while preserving their grading. We say that $A$ and $B$ are \emph{$(H,K)$-quantum-symmetrically equivalent} (denoted $\HKQS$ equivalent) if there is a monoidal equivalence \[F: \comod(H)\to \comod(K)\]
   such that $F$ sends $A$ to $B$ as comodule algebras. 
\end{defn}

\begin{remark}
    By our definition, two connected graded, locally finite algebras $A,B$ are QS equivalent if and only if they are $(\underline{\rm aut}^r(A),\underline{\rm aut}^r(B))$-QS equivalent for their associated Manin's universal quantum groups. 
\end{remark}

We now show that for any two connected graded algebras that are locally finite,  $\HKQS$ equivalence implies QS equivalence. We point out that the assumption of local finiteness in the result below is necessary only to ensure the existence of Manin's universal quantum groups.

\begin{thm}
\label{RQS}
Let $A$ and $B$ be two connected graded algebras that are locally finite. If $A$ and $B$ are $\HKQS$ equivalent, for some Hopf algebras $H$ and $K$ right coacting on $A$ and $B$ respectively, preserving their grading, then $A$ and $B$ are QS equivalent.   
\end{thm}

\begin{proof}
Suppose $A$ and $B$ are $\HKQS$ equivalent for some Hopf algebras $H$ and $K$ right coacting on $A$ and $B$ respectively. Let $T$ be the $H$-$K$-bi-Galois object that induces the monoidal equivalence 
\[F(-)=(-)\cotens_HT: \comod(H)\to \comod(K),\]
such that $F(A)=B$ as comodule algebras. By the universal property of $\underline{\rm aut}^r(A)$, there is a unique Hopf algebra map $f: \underline{\rm aut}^r(A)\to H$ making the corresponding diagram of coactions on $A$ commute. Thus, we have an induced monoidal functor $f_*: \comod(\underline{\rm aut}^r(A))\to \comod(H)$ that is $\kk$-linear, exact and faithful. Consider the left $\underline{\rm aut}^r(A)$-comodule algebra $P \coloneqq \underline{\rm aut}^r(A)\cotens_HT$, where we view $\underline{\rm aut}^r(A)$ as a right $H$-comodule via $f_*$. It is clear that for any right $\underline{\rm aut}^r(A)$-comodule $M$, 
\[M\cotens_{\underline{\rm aut}^r(A)} P\cong (M\cotens_{\underline{\rm aut}^r(A)} {\underline{\rm aut}^r(A)}) \cotens_H T\cong f_*(M)\cotens_HT,\] 
which implies that cotensoring with $P$ is a fiber functor from $\comod(\underline{\rm aut}^r(A))$ to ${\rm Vec}(\kk)$. By Ulbrich’s theorem \cite{Ulbrich} (cf. \cite[Theorem 1.14]{Bichon2014}), $P$ is a left $\underline{\rm aut}^r(A)$-Galois object. Moreover, by \Cref{Schauenburg}, there is a unique Hopf algebra $L$ where $P$ becomes an $\underline{\rm aut}^r(A)$-$L$-bi-Galois object such that cotensoring with $P$ induces a monoidal equivalence between $\comod(\underline{\rm aut}^r(A))$ and $\comod(L)$. Additionally, we have
\[
A\cotens_{\underline{\rm aut}^r(A)}P=A\cotens_{\underline{\rm aut}^r(A)}(\underline{\rm aut}^r(A)\cotens_HT)\cong (A\cotens_{\underline{\rm aut}^r(A)}\underline{\rm aut}^r(A))\cotens_HT\cong f_*(A)\cotens_HT=A\cotens_HT=B,
\] 
where the right coaction of $L$ on $B$ is induced from the right coaction of $L$ on $P$. Therefore, cotensoring with $P$ gives a $(\underline{\rm aut}^r(A),L)$-QS equivalence between $A$ and $B$. 

It remains to show that $L\cong \underline{\rm aut}^r(B)$ as Hopf algebras. Without loss of generality, we can take $H=\underline{\rm aut}^r(A)$ and show that $K\cong \underline{\rm aut}^r(B)$. Let $S$ be the $K$-$H$-bi-Galois object such that $(-)\cotens_KS:\comod(K)\to \comod(H)$ gives the monoidal inverse of $(-)\cotens_HT: \comod(H)\to \comod(K)$ while sending $B$ to $A$ as comodule algebras. As noted above, there is a unique Hopf algebra map $g: \underline{\rm aut}^r(B)\to K$ where the corresponding coactions on $B$ factor through. Moreover, by a similar argument, there is a unique Hopf algebra $X$ such that cotensoring with the $\underline{\rm aut}^r(B)$-$X$-bi-Galois object $Q \coloneqq \underline{\rm aut}^r(B)\cotens_KS$ yields an $(\underline{\rm aut}^r(B),X)$-QS equivalence between $B$ and $A$. Observe that since $H$ right coacts on $S$, $H$ also right coacts on $Q$. By the universal property of the left $\underline{\rm aut}^r(B)$-Galois object $Q$ (cf. \Cref{Schauenburg}), we have a unique Hopf algebra map $g': X\to H$ such that the $H$-coaction on $Q$ factors through that of $X$ via $g'$. We claim  the diagram of monoidal functors
\[
\xymatrix{
\comod(\underline{\rm aut}^r(B)\ar[d]_-{g_*})\ar[rr]^-{(-)\cotens_{\underline{\rm aut}^r(B)}Q} && \comod(X)\ar[d]^-{g'_*}\\
\comod(K)\ar[rr]^-{(-)\cotens_KS}&& \comod(H)
}
\]
naturally commutes. Indeed, any right $\underline{\rm aut}^r(B)$-comodule $V$ can be viewed as a right $K$-comodule via $g: \underline{\rm aut}^r(B)\to K$. Moreover, we get isomorphisms of vector spaces  
\[
g_*'(V\cotens_{\underline{\rm aut}^r(B)}Q) \cong V\cotens_{\underline{\rm aut}^r(B)}(\underline{\rm aut}^r(B)\cotens_KS) 
\cong (V\cotens_{\underline{\rm aut}^r(B)}\underline{\rm aut}^r(B))\cotens_KS\cong g_*(V)\cotens_KS,
\]
where the corresponding $H$-coaction is induced by the $H$-coaction on $S$. It is clear that these isomorphisms are natural in terms of the $\underline{\rm aut}^r(B)$-comodule $V$, that is, for any other right $\underline{\rm aut}^r(B)$-comodule $W$ and $\underline{\rm aut}^r(B)$-comodule map $f: V\to W$, we have the following commutative diagram   
\[
\xymatrix{
g_*'(V\cotens_{\underline{\rm aut}^r(B)}Q)\ar[d]_-{g_*'(f \cotens_{\underline{\rm aut}^r(B)}Q)}\ar[r]^-{\cong} & g_*(V)\cotens_KS\ar[d]^-{g_*(f)\cotens_KS}\\
g_*'(W\cotens_{\underline{\rm aut}^r(B)}Q)\ar[r]^-{\cong} & g_*(W)\cotens_KS\, .
}
\]
This proves our claim. In particular, $B\cotens_{\underline{\rm aut}^r(B)}Q\cong B\cotens_KS=A$, which implies that $A$ is an $X$-comodule algebra. By our assumption $H=\underline{\rm aut}^r(A)$, so by the universal property of $\underline{\rm aut}^r(A)$, there is a unique Hopf algebra map $h: H\to X$ where $g'\circ h=\id$ on $H$. 

Let $R \coloneqq X\cotens_HT$, which is a left $X$-Galois object since $T$ is a left $H$-Galois object. One can check that $(-)\cotens_XR: \comod(X)\to \comod(\underline{\rm aut}^r(B))$ is the monoidal inverse of $(-)\cotens_{\underline{\rm aut}^r(B)} Q: \comod(\underline{\rm aut}^r(B))\to \comod(X)$ since 
 \[
 (V\cotens_X R)\cotens_{\underline{\rm aut}^r(B)} Q\cong (g_*'(V)\cotens_HT)\cotens_KS\cong g_*'(V)\]
 and 
 \[(W\cotens_{\underline{\rm aut}^r(B)} Q)\cotens_XR\cong (g_*(W)\cotens_KS)\cotens_HT \cong g_*(W)
 \]
 for any $V\in \comod(X)$ and $W\in \comod(\underline{\rm aut}^r(B))$. Because $g'\circ h=\id$ on $H$, it is clear that $h_*(H)\cotens_XR=g'_*\circ h_*(H)\cotens_{H} T=(g'\circ h)_*(H)\cotens_HT=T\in \comod(\underline{\rm aut}^r(B))$. By the universal property of $T$ as a right $K$-Galois object, we get a unique Hopf algebra map $h': K\to \underline{\rm aut}^r(B)$. In summary, we have the following commutative diagram of monoidal functors:
\[
\xymatrix{
\comod(H)\ar[rr]^-{(-)\cotens_HT}\ar[d]_-{h_*} & &\comod(K)\ar[d]^-{h'_*}\\
\comod(X)\ar[rr]^-{(-)\cotens_XR}\ar[d]_-{g'_*} & &\comod(\underline{\rm aut}^r(B))\ar[d]^-{g_*}\\
\comod(H)\ar[rr]^-{(-)\cotens_HT} &&\comod(K).
}
\]
Since $g'_*\circ h_*=(g'\circ h)_*=\id_*$, we get $g_*\circ h'_*=(g\circ h')_*$ is an auto-equivalence on $\comod(K)$ and hence on $\comod_{\fd}(K)$, which shows that $g\circ h'$ is an isomorphism on $K$ by the Tannaka--Krein reconstruction theorem, and hence $g$ is surjective. By the universal coaction of $\underline{\rm aut}^r(B)$ on $B$, we also have $h'\circ g=\id$ on $\underline{\rm aut}^r(B)$ and then $g$ is injective. Hence, $g:\underline{\rm aut}^r(B)\to K$ is an isomorphism of Hopf algebras. This shows that $A$ and $B$ are QS equivalent.
\end{proof}

We recall some properties of connected graded algebras that are invariant under $\HKQS$ equivalence.

\begin{thm}[{\cite[Corollary 3.2.7 and Theorem 3.2.10]{HNUVVW2024-2}}]
\label{PQS}
Let $A,B$ be two connected graded algebras, finitely generated in degree one. Suppose $A,B$ are $\HKQS$ equivalent. We have:
 \begin{enumerate}
     \item[(1)] ${\rm gl.dim}(A)={\rm gl.dim}(B)$;
     \item[(2)] $A$ is $N$-Koszul if and only if $B$ is $N$-Koszul.
 \end{enumerate}
Further, assume $A$ is noetherian and $H$ has a bijective antipode. If $A$ is AS-regular of dimension $d$, then so is $B$.
\end{thm}

\subsection{$\mathcal{GL}$-type and $\mathcal{SL}$-type quantum-symmetric equivalences}

For any integer $m \ge 2$, we investigate the cogroupoid $\mathcal{GL}_m$ constructed in \Cref{cogp}. Let $e\in U^{\otimes m},f\in V^{\otimes m}$ be two nondegenerate twisted superpotentials where $U, V$ are two $\kk$-vector spaces with fixed bases $\{u_1, \ldots, u_p\}$ and $\{v_1, \ldots, v_q\}$, respectively. We also write $\{u^1, \ldots, u^p\}$ and $\{v^1, \ldots, v^q\}$ for the corresponding dual bases of $U^*$ and $V^*$. Suppose the bi-Galois object $\mathcal{GL}_m(e,f)\neq 0$. Then, a monoidal equivalence between two right comodule categories:
\begin{align}
\label{Fcogroupoid}
(F,\xi): \comod(\mathcal{GL}_m(e))~&\to~ \comod(\mathcal{GL}_m(f)) \\
W~&\mapsto~ F(W)\coloneqq  W \cotens_{\mathcal{GL}_m(e)}\,\mathcal{GL}_m(e,f), \notag
\end{align}
can be given by the cotensor product with the bi-Galois object $\mathcal{GL}_m(e,f)$, and $\xi$ is a natural isomorphism (as in \cite[Proposition 1.4]{Ulbrich1987} for left comodules) 
\begin{align*}
\xi_{W,M}: (W\cotens_{\mathcal{GL}_m(e)}\mathcal{GL}_m(e,f))\otimes (M\cotens_{\mathcal{GL}_m(e)}\mathcal{GL}_m(e,f))~&\cong~ (W\otimes M)\cotens_{\mathcal{GL}_m(e)}\mathcal{GL}_m(e,f)\\
\left(\sum w_i\otimes a_i\right)\otimes \left(\sum m_j\otimes b_j\right)~&\mapsto~ \sum w_i\otimes m_j\otimes a_ib_j
\end{align*}
satisfying certain compatible conditions (cf. \cite[\S 2.4]{EGNO2015}) for $w_i \in W, m_j \in M$ and $a_i,b_j \in \mathcal{GL}_m(e,f)$. In our context, the cotensor product is explicitly given by 
\begin{align}\label{coten}
W\cotens_{\mathcal{GL}_m(e)}\mathcal{GL}_m(e,f) \coloneqq \left\{x\in W\otimes \mathcal{GL}_m(e,f)\,|\, (\rho_W\otimes \id_{\mathcal{GL}_m(e,f)})(x)=(\id_W \otimes \Delta_{e,f}^e)(x)\right\},
\end{align}
where $\rho_W: W\to W\otimes \mathcal{GL}_m(e)$ is the right $\mathcal{GL}_m(e)$-comodule structure map on $W$. Moreover, the right $\mathcal{GL}_m(f)$-comodule structure on $F(W)= W\cotens_{\mathcal{GL}_m(e)}\mathcal{GL}_m(e,f)$ comes from that of $\mathcal{GL}_m(f)$ on $\mathcal{GL}_m(e,f)$ via $\Delta_{e,f}^f$.

We denote by $U_e\coloneqq U$ the right comodule over $\mathcal{GL}_m(e)$ given via the $\mathcal{GL}_m(e)$-coaction $u_i\mapsto \sum_{1\leq j\leq p} u_j\otimes a_{ji}^{e,e} \in U_e \otimes \mathcal{GL}_m(e)$. Similarly, we denote by $V_f\coloneqq V$ the analogous right comodule over $\mathcal{GL}_m(f)$. Using \eqref{eq:alg}, it is straightforward to check that 
\begin{align*}
  e: \kk(D^{e,e})^{-1} &\to U_e^{\otimes m}\\
 (D^{e,e})^{-1} &\mapsto \sum_{1\leq i_1,\ldots,i_m\leq p}e_{i_1\cdots i_m}u_{i_1}\otimes \cdots \otimes u_{i_m}
\end{align*}
is a right $\mathcal{GL}_m(e)$-comodule map, where $\kk(D^{e,e})^{-1}$ is the one-dimensional right $\mathcal{GL}_m(e)$-comodule given by the group-like element $(D^{e,e})^{-1}$. The $\kk$-vector space dual $U_e^*$ equipped with the right $\mathcal{GL}_m(e)$-comodule structure given by $u^i\mapsto \sum_{1\leq j\leq p}u^j\otimes b_{ij}$ is the left dual of $U_e$ in the category $\comod(\mathcal{GL}_m(e))$. Moreover, the dual map 
\[
e^*~=~\left((U_e^*)^{\otimes m}\xrightarrow{\id \otimes {\rm coev} } (U_e^*)^{\otimes m} \otimes \kk (D^{e,e})^{-1} \otimes \kk (D^{e,e})\xrightarrow{\id \otimes e \otimes \id } (U_e^*)^{\otimes m} \otimes U_e^{\otimes m} \otimes \kk (D^{e,e})\xrightarrow{{\rm ev}\otimes \id }\kk D^{e,e}\right)
\]
is again $\mathcal{GL}_m(e)$-collinear and is given by
\begin{align*}
    e^*(u^{i_1}\otimes u^{i_m})&\,=\left\langle u^{i_1}\otimes \cdots \otimes u^{i_m}, \sum\limits_{1\leq j_1,\ldots,j_m\leq p} e_{j_1\cdots j_m}u_{j_1}\otimes \cdots \otimes u_{j_m}\right \rangle\, D^{e,e}\\
    &\,= \sum\limits_{1\leq j_1,\ldots,j_m\leq p}\langle u^{i_1},u_{j_m}\rangle\cdots \langle u^{i_m},u_{j_1}\rangle\,e_{j_1\cdots j_m}\, D^{e,e}\\
    &\,= e_{i_m\cdots i_1}\, D^{e,e}.
\end{align*}
Similar formulas hold for the $\mathcal{GL}_m(f)$-collinear maps $f: \kk(D^{f,f})^{-1}\to V_f^{\otimes m}$ and $f^*: (V_f^*)^{\otimes m}\to \kk D^{f,f}$.

In what follows, when we say that the morphism $ e: \kk(D^{e,e})^{-1} \to U_e^{\otimes m}$ is nondegenerate, we are viewing $\Imm(e)$ as a nondegenerate twisted superpotential as in \Cref{defn:superpotential}(1)(ii).

\begin{lemma}
\label{lem:comoduleGL}
Let $2\leq N\leq m$ be integers,  $e\in U^{\otimes m}$ a nondegenerate twisted superpotential, and retain the above notation. Then the corresponding superpotential algebra $A(e,N)=TU_e/\langle\partial^{m-N}(\kk e)\rangle$ is a right comodule algebra over $\mathcal{GL}_m(e)$ and over $\mathcal{SL}_m(e)$.
\end{lemma}
\begin{proof}
  Since $U_e$ is a right $\mathcal{GL}_m(e)$-comodule, it suffices to show that $\partial^{m-N}(\kk e)\subset U_e^{\otimes N}$ is a subcomodule over $\mathcal{GL}_m(e)$.
  Observe that for any integer $n \geq 1$ and any $\mathcal{GL}_m(e)$-subcomodule $\iota: W\subset U_e^{\otimes n}$, it is straightforward to check that the subspace $\partial(W) \subset U_e^{\otimes (n-1)}$ can be given by  
  \[
  \partial(W)={\rm Image}\left(U_e^*\otimes W\xrightarrow{\id \otimes \iota} U_e^*\otimes U_e^{\otimes n}\xrightarrow{\cong } (U_e^*\otimes U_e)\otimes U_e^{\otimes (n-1)}\xrightarrow{{\rm ev}\otimes \id} U_e^{\otimes (n-1)}\right).
  \]
  Since $\partial(W)$ is the image of the composition of all right $\mathcal{GL}_m(e)$-comodule maps, we have $\partial(W)$ is a $\mathcal{GL}_m(e)$-subcomodule of $U_e^{\otimes (n-1)}$. Iteratively, $\partial^i(W)=\partial(\partial^{i-1}(W))\subset U_e^{\otimes (n-i)}$ are subcomodules for all $2\leq i\leq n$. Applying this observation to the $\mathcal{GL}_m(e)$-comodule embedding $e: \kk(D^{e,e})^{-1} \to U_e^{\otimes m}$, it implies that $\partial^{m-N}(\kk e)=\partial^{m-N}(\kk (D^{e,e})^{-1})$ is a $\mathcal{GL}_m(e)$-subcomodule of $U_e^{\otimes N}$. Finally, the $\mathcal{SL}_m(e)$-coaction on $A(e,N)$ is induced from the Hopf quotient map $\pi:\mathcal{GL}_m(e)\to \mathcal{SL}_m(e)=\mathcal{GL}_m(e)/(D^{e,e}-1)$.
\end{proof}

In light of \Cref{def:HKQS}, we propose the following notation for the two types of quantum-symmetric equivalence for superpotential algebras regarding Hopf coactions associated with $\mathcal{GL}_m$ and $\mathcal{SL}_m$.

\begin{notation}
\label{GLQSandSLQS}
Let $2\leq N\leq m$ be integers. Suppose $A(e,N)$ and $A(f,N)$ are two superpotential algebras given by nondegenerate twisted superpotentials $e\in U^{\otimes m}$ and $f\in V^{\otimes m}$, respectively. 
\begin{itemize}
      \item[(a)] We say $A(e,N)$ and $A(f,N)$ are \emph{$\mathcal{GL}$ quantum-symmetrically equivalent} (denoted $\GLQS$ equivalent) if they are $(\mathcal{GL}_m(e),\mathcal{GL}_m(f))$-QS equivalent.
       \item[(b)] We say $A(e,N)$ and $A(f,N)$ are \emph{$\mathcal{SL}$ quantum-symmetrically equivalent} (denoted $\SLQS$ equivalent) if they are $(\mathcal{SL}_m(e),\mathcal{SL}_m(f))$-QS equivalent.
\end{itemize}
\end{notation}

We now provide a condition for two superpotential algebras to be $\GLQS$ equivalent.

\begin{thm}
\label{SQSMTE}
Let $2\leq N\leq m$ be integers and $e\in U^{\otimes m}$ and $f\in V^{\otimes m}$ be two nondegenerate twisted superpotentials. If $\mathcal{GL}_m(e,f)\neq 0$, then the corresponding superpotential algebras $A(e,N)$ and $A(f,N)$ are $\GLQS$ equivalent. 
\end{thm}

\begin{proof}
Retain the above notation, with $\{u_1, \ldots, u_p\}$ and $\{v_1, \ldots, v_q\}$ being bases for the $\kk$-vector spaces $U$ and $V$, respectively. To show that $A(e,N)$ and $A(f,N)$ are $\GLQS$ equivalent, it is enough to show that $F(A(e,N))=A(f,N)$, where $F$ is the cotensor product given in \eqref{Fcogroupoid}. By the construction of the superpotential algebras, it is equivalent to show that  $F(U_e)=V_f$, and the morphism $e: \kk (D^{e,e})^{-1}\to U_e^{\otimes m}$ is sent to the morphism $f: \kk (D^{f,f})^{-1}\to V_f^{\otimes m}$. 

We first follow the argument in the proof of \cite[Corollary 3.5]{Bichon2001} to show that $F(U_e)=V_f$. As a consequence of \cite[Theorem 2.12]{Bichon2001}, $\mathcal{GL}_m(f,e)\neq 0$ and we have the inverse monoidal equivalence 
\begin{align*}
(G,\zeta): \comod(\mathcal{GL}_m(f))~&\to~ \comod(\mathcal{GL}_m(e)) \\
M~&\mapsto~ G(M)\coloneqq  M\,\cotens_{\mathcal{GL}_m(f)}\,\mathcal{GL}_m(f,e).
\end{align*}
It is straightforward to check that the linear map
\begin{align*}
  \nu_f: V_f\to&\, F(U_e)=U_e\cotens_{\mathcal{GL}_m(e))}\mathcal{GL}(e,f)\\
  v_j\mapsto&\,  \sum_{1\leq i\leq p} u_i\otimes a_{ij}^{e,f}, \quad \text{for all $1 \leq j \leq q$,}
\end{align*}
is $\mathcal{GL}_m(f)$-collinear. Similarly, we have a $\mathcal{GL}_m(e)$-collinear map $\nu_e: U_e\to G(V_f)$. We then obtain a sequence of collinear maps
\[
V_f\xrightarrow{\nu_f} F(U_e)\xrightarrow{F(\nu_e)} FG(V_f)\xrightarrow{\cong }V_f
\]
whose composition is the identity map, as shown by the following concrete computation: 
\begin{align*}
 & V_f\xrightarrow{} U_e\cotens_{\mathcal{GL}_m(e)}\mathcal{GL}_m(e,f)\xrightarrow{} \left(V_f\cotens_{\mathcal{GL}_m(f)}\mathcal{GL}_m(f,e)\right)\cotens_{\mathcal{GL}_m(e)}\mathcal{GL}_m(e,f) \\ 
 & \quad \xrightarrow{\cong} V_f\cotens_{\mathcal{GL}_m(f)}(\mathcal{GL}_m(f,e)\cotens_{\mathcal{GL}_m(e)}\mathcal{GL}_m(e,f))
\xrightarrow{\cong} V_f\cotens_{\mathcal{GL}_m(f)}\mathcal{GL}_m(f)\xrightarrow{\cong }V_f, \text{ via}  \\ 
& v_j\mapsto \sum_{1\leq i\leq p} u_i\otimes a_{ij}^{e,f}\mapsto \sum_{1\leq i\leq p} \sum_{1\leq k\leq q} v_k\otimes a_{ki}^{f,e}\otimes a_{ij}^{e,f}=\sum_{1\leq k\leq q}v_k\otimes \Delta_{f,f}^e(a_{kj})\mapsto v_j.
\end{align*}
Hence $\nu_f$ is injective and $F(\nu_e)$ is surjective, and $\nu_f$ is surjective by a symmetric argument. So, we can conclude $F(U_e)=V_f$.

Next, to show that $\xi^{-1}F(e)=f$, we note that 
\[F(\kk (D^{e,e})^{-1})=\kk (D^{e,e})^{-1}\cotens_{\mathcal{GL}_m(e)}\mathcal{GL}_m(e,f)=\kk (D^{e,e})^{-1}\otimes (D^{e,f})^{-1}\cong \kk (D^{f,f})^{-1}.\]
We claim the diagram 
\[
\xymatrix{
\kk (D^{f,f})^{-1}\ar[r]^-{f}\ar[d]_-{\cong} & V_f^{\otimes m}\ar[r]^-{=} & F(U_e)^{\otimes m}\ar[r]^-{\cong } & (U_e\cotens_{\mathcal{GL}_m(e)}\mathcal{GL}_m(e,f))^{\otimes m}\ar[d]^-{\xi} \\
 \kk (D^{e,e})^{-1} \cotens_{\mathcal{GL}_m(e)}\mathcal{GL}_m(e,f)\ar[rrr]^-{e\otimes \id} & && (U_e^{\otimes m})\cotens_{\mathcal{GL}_m(e)}\mathcal{GL}_m(e,f) 
}
\]
commutes. By tracing the diagram at the level of elements, we obtain the composition of maps on the top row and right column as 
\begin{align*}
    (D^{f,f})^{-1}&\overset{f}{\mapsto} \sum_{1\leq i_1,\ldots,i_m\leq q} f_{i_1\cdots i_m}v_{i_1}\otimes \cdots \otimes v_{i_m}\overset{\cong }{\mapsto} \sum_{\stackrel{1\leq i_1,\ldots,i_m\leq q}{1\leq j_1,\ldots,j_m\leq p}} f_{i_1\cdots i_m}(u_{j_1}\otimes a_{j_1i_1}^{e,f})\otimes \cdots \otimes (u_{j_m}\otimes a_{j_mi_m}^{e,f})\\
    &\overset{\xi}{\mapsto} \sum_{1\leq j_1,\ldots,j_m\leq p} u_{j_1}\otimes \cdots \otimes u_{j_m}\otimes \Big( \sum_{1\leq i_1,\ldots,i_m\leq q}f_{i_1\cdots i_m}a_{j_1i_1}^{e,f}\cdots a_{j_mi_m}^{e,f} \Big) \\
\overset{\eqref{eq:alg}}&{=}  \sum_{1\leq j_1,\ldots,j_m\leq p} e_{j_1\cdots j_m}u_{j_1}\otimes \cdots \otimes u_{j_m}\otimes (D^{e,f})^{-1}.
\end{align*}
For the composition of maps on the left column and bottom row, we have 
\begin{align*}
    \kk (D^{f,f})^{-1}&\overset{\cong }{\mapsto} (D^{e,e})^{-1}\otimes (D^{e,f})^{-1}\overset{e\otimes \id }{\mapsto}  \sum_{1\leq j_1,\ldots,j_m\leq p}e_{j_1\cdots j_m} u_{j_1}\otimes \cdots \otimes u_{j_m} \otimes (D^{e,f})^{-1}.    
\end{align*}
They are equal, which shows the commutativity of the above diagram. Moreover, it is straightforward to see that the diagram 
\[
\xymatrix{
\kk (D^{f,f})^{-1}\cong F(\kk(D^{e,e})^{-1})\ar[r]^-{F(e)}\ar[d]_-{\cong} & F(U_e^{\otimes m})\ar[r]^-{\xi^{-1}}\ar[d]_-{\cong}& F(U_e)^{\otimes m}\ar[r]^-{\cong}\ar[d]_-{\cong } & V_f^{\otimes m}\ar[dl]_-{\cong}\\
 \kk (D^{e,e})^{-1} \cotens_{\mathcal{GL}_m(e)}\mathcal{GL}_m(e,f)\ar[r]^-{e\otimes \id} & (U_e^{\otimes m})\cotens_{\mathcal{GL}_m(e)}\mathcal{GL}_m(e,f)\ar[r]^-{\xi^{-1}}& (U_e\cotens_{\mathcal{GL}_m(e)}\mathcal{GL}_m(e,f))^{\otimes m}& 
}
\]
commutes. Thus, we get $\xi^{-1}F(e)=\xi^{-1}(e\otimes \id)=f$. 
\end{proof}

\subsection{The Artin--Schelter regular component of the cogroupoid $\mc{GL}_m$}

We propose a concept of traceable superpotentials, which allows the converse of \Cref{SQSMTE} to also hold.

\begin{defn}
\label{defn:traceable}
A nondegenerate twisted superpotential $f\in V^{\otimes m}$ is said to be \emph{$L$-traceable} for some integers $2\leq L\leq m$ if 
\[
\kk f=\bigcap_{i+j=m-L} V^{\otimes i}\otimes \partial^{m-L}(\kk f)\otimes V^{\otimes j}\subset V^{\otimes m}.
\] 
\end{defn}

It is known, by \cite[Proposition 2.12]{MS2016}, that superpotentials for $N$-Koszul AS-regular algebras are always $N$-traceable. It is an open question whether the converse also holds.

\begin{proposition}
\label{traceable}
Let $2 \leq N\leq m$ be integers and $e\in U^{\otimes m}$ and $f\in V^{\otimes m}$ be two nondegenerate twisted superpotentials.
If $A(e,N)$ and $A(f,N)$ are $\GLQS$ equivalent, then $A(e,N')$ and $A(f,N')$ are $\GLQS$ equivalent for any $2\leq N'\leq N$. If we additionally assume that $e,f$ are $L$-traceable for some $2\leq L\leq N$, then $A(e,N')$ and $A(f,N')$ are $\GLQS$ equivalent for any $2\leq N'\leq m$. 
\end{proposition}

\begin{proof}
By definition, there is a monoidal equivalence
\begin{align*}
(F,\xi): \comod(\mathcal{GL}_m(e))\to \comod(\mathcal{GL}_m(f)),
\end{align*}
where $F$ sends $A(e,N)=TU/\langle\partial^{m-N}(\kk e)\rangle$ to $A(f,N)=TV/\langle\partial^{m-N}(\kk f)\rangle$ as comodule algebras. Therefore, we can assume that $F(U)=V$ and $\xi^{-1} F(\partial^{m-N}(\kk e))=\partial^{m-N}(\kk f)\subseteq V^{\otimes N}$.  

For any integer $n \geq 1$ let $\iota: W\subset U^{\otimes n}$ be a subcomodule over $\mathcal{GL}_m(e)$. We claim that $\xi^{-1} F(\partial W)=\partial(\xi^{-1} F(W))\subseteq V^{\otimes (n-1)}$. This follows from the commutativity of the diagram
\[\small
\xymatrix{
F(U^*\otimes W)\ar[r]^-{F(\id \otimes \iota)}\ar[d]_-{\xi^{-1}} & F(U^*\otimes U^{\otimes n})\ar[d]_-{\xi^{-1}}\ar[r]^-{F(\cong)}& F((U^*\otimes U)\otimes U^{n-1})\ar[r]^-{F({\rm ev}\otimes \id)}\ar[d]_-{\xi^{-1}}& F(U^{\otimes (n-1)})\ar@{=}[d]\\
F(U^*)\otimes F(W)\ar[d]_-{\id \otimes F(\iota)}\ar[r]^-{\id \otimes F(\iota)} &F(U^*)\otimes F(U^{\otimes n})\ar[d]_-{\id \otimes \xi^{-1}}&F(U^*\otimes U)\otimes F(U^{\otimes (n-1)})\ar[d]_-{\xi^{-1}\otimes \xi^{-1}}\ar[r]^-{F({\rm ev})\otimes \id}& F(U^{\otimes (n-1)})\ar@{=}[d]\\
F(U^*)\otimes F(U^{\otimes n})\ar[d]_-{\widetilde{\xi}\otimes \id}\ar[r]^-{\id \otimes \xi^{-1}} &F(U^*)\otimes F(U)^{\otimes n}\ar[d]_-{\widetilde{\xi}\otimes \id}\ar[r]^-{\cong}& (F(U^*)\otimes F(U))\otimes F(U)^{\otimes (n-1)}\ar[d]_-{\widetilde{\xi} \otimes \id\otimes \id}&F(U^{\otimes (n-1)})\ar[d]^-{\xi^{-1}}\\
F(U)^*\otimes F(U^{\otimes n})\ar[r]^-{\id \otimes \xi^{-1}} &F(U)^*\otimes F(U)^{\otimes n}\ar[r]^-{\cong}& (F(U)^*\otimes F(U))\otimes F(U)^{\otimes (n-1)}\ar[r]^-{{\rm ev}\otimes \id}& F(U)^{\otimes (n-1)}
}
\]
where $\widetilde{\xi}$ is defined in \eqref{xitilde}. Observe that the compositions of maps appearing on the outer rectangle show that $\xi^{-1} F(\partial W)$ and $\partial(\xi^{-1} F(W))$ share the same image in $V^{\otimes (n-1)}$. By iteration, one can show that for any $2 \leq N' \leq N$,
\begin{align*}
  \xi^{-1}F(\partial^{m-N'}(\kk e))&\,=\xi^{-1}F(\partial^{N-N'}(\partial^{m-N}(\kk e)))=\partial^{N-N'}(\xi^{-1} F(\partial^{m-N}(\kk e)))\\
  &\,=\partial^{N-N'}(\partial^{m-N}(\kk f))=\partial^{m-N'}(\kk f).
 \end{align*}

Consequently, $F$ sends $A(e,N')=TU/\langle\partial^{m-N'}(\kk e)\rangle$ to $A(f,N')=TV/\langle\partial^{m-N'}(\kk f)\rangle$ as comodule algebras. Therefore, $A(e,N')$ and $A(f,N')$ are $\GLQS$ equivalent for any $2 \leq N' \leq N \leq m$.

For the second statement of the lemma, now we additionally assume $e,f$ are $L$-traceable for some $2\leq L\leq N$. Thus, we have
\begin{align*}
    \xi^{-1} F(\kk e)&\,=\xi^{-1} F\left(\bigcap_{i+j=m-L} U^{\otimes i}\otimes \partial^{m-L}(\kk e)\otimes  U^{\otimes j} \right)=\bigcap_{i+j=m-L} V^{\otimes i}\otimes \xi^{-1} F(\partial^{m-L}(\kk e))\otimes V^{\otimes j}\\
    &\, =\bigcap_{i+j=m-L} V^{\otimes i}\otimes \partial^{m-L}(\kk f)\otimes V^{\otimes j}=\kk f.
\end{align*}
This implies that $A(e,m)=TU/(\kk e)$ and $A(f,m)=TV/(\kk f)$ are $\GLQS$ equivalent. The result follows immediately from our previous discussion. 
\end{proof}

Now assume nondegenerate twisted superpotentials $e\in U^{\otimes m},f\in V^{\otimes m}$ are $L$-traceable for some $2\leq L\leq N\leq m$ and $A(e,N)$ and $A(f,N)$ are $\GLQS$ equivalent via some monoidal equivalence
\begin{align}
\label{MEUQ}
(F,\xi): \comod(\mathcal{GL}_m(e))\to \comod(\mathcal{GL}_m(f)),
\end{align}
where $F$ sends $A(e,N)$ to $A(f,N)$ as comodule algebras. By \Cref{traceable}, $\partial^{m-L}(\kk e)\subseteq U_e^{\otimes L}$ and $\partial^{m-L}(\kk f)\subseteq V_f^{\otimes L}$ are subcomodules over $\mathcal{GL}_m(e)$ and $\mathcal{GL}_m(f)$, respectively.  As a consequence of $e,f$ being $L$-traceable, we have 
\begin{align*}
 \kk (D^{e,e})^{-1}~&\overset{e}{=}~\bigcap_{i+j=m-L} U_e^{\otimes i}\otimes \partial^{m-L}(\kk e)\otimes U_e^{\otimes j}\subset U_e^{\otimes m} \quad \text{ and}\\
\kk (D^{f,f})^{-1}~&\overset{f}{=}~\bigcap_{i+j=m-L} V_f^{\otimes i}\otimes \partial^{m-L}(\kk f)\otimes V_f^{\otimes j}\subset V_f^{\otimes m},   
\end{align*}
are one-dimensional invertible comodules over $\mathcal{GL}_m(e)$ and $\mathcal{GL}_m(f)$, respectively. Moreover, we obtain $\xi^{-1} F(\kk e)=\kk f$ such that 
\begin{align}\label{Tf}
    f~=~\left(\kk (D^{f,f})^{-1}=F(\kk(D^{e,e})^{-1})\xrightarrow{F(e)}F(U_e^{\otimes m})\xrightarrow{\xi^{-1}} F(U_e)^{\otimes m}\xrightarrow{=}V_f^{\otimes m}\right).
\end{align}

Under this setting, we now show the converse of \Cref{SQSMTE} also holds.

\begin{thm}
\label{SQSBG}
If $A(e,N)$ and $A(f,N)$ are $\GLQS$ equivalent for some $L$-traceable nondegenerate twisted superpotentials $e,f$ with integers $2\leq L\leq N\leq m$, then $\mathcal{GL}_m(e,f)\neq 0$. Moreover, any monoidal equivalence from $\comod(\mathcal{GL}_m(e))$ to $\comod(\mathcal{GL}_m(f))$ that sends $A(e,N)$ to $A(f,N)$ as comodule algebras is naturally isomorphic to $(-)\cotens_{\mathcal{GL}_m(e)}\mathcal{GL}_m(e,f)$. 
\end{thm}

\begin{proof}
    Since any monoidal equivalence $F:\comod(\mathcal{GL}_m(e))\to \comod(\mathcal{GL}_m(f))$ in \eqref{MEUQ} is given by cotensoring with some nonzero $\mathcal{GL}_m(e)$-$\mathcal{GL}_m(f)$-bi-Galois object $T$, it suffices to show that there exists an $\mathcal{GL}_m(e)$-$\mathcal{GL}_m(f)$-bicomodule algebra isomorphism $\pi: \mathcal{GL}_m(e,f)\cong T$. 

Suppose $e \in U^{\otimes m}, f \in V^{\otimes m}$ where $\dim U=p$ and $\dim V=q$. As usual, write bases $\{u_1,\ldots,u_p\}$ and $\{v_1,\ldots,v_q\}$ for $U_e$ and $V_f$ with dual bases $\{u^1,\ldots,u^p\}$ and $\{v^1,\ldots,v^q\}$ for $U_e^*$ and $V_f^*$, respectively. Since $F$ sends $A(e,N)=TU_e/\langle\partial^{m-N}(\kk e)\rangle$ to $A(f,N)=TV_f/\langle\partial^{m-N}(\kk f)\rangle$ as comodule algebras, we have $F(U_e)=U_e\cotens_{\mathcal{SL}_m(e)}T=V_f$.  Then there are elements $(y_{ji})_{1\leq  j\leq p,1\leq i\leq q}$ of $T$ such that $v_i=\sum_{1\leq j\leq p} u_j\otimes y_{ji}$ for all $1 \leq i \leq q$. We have $U_e^*\cotens_{\mathcal{SL}_m(e)}T= F(U_e^*)\cong F(U_e)^*=V_f^*$ where the middle isomorphism is given by the duality transformation $\widetilde{\xi}$. Therefore, there are elements $(z_{ij})_{1\leq i\leq q,1\leq j\leq p}$ such that $v^i=\sum_{1\leq j\leq p} u^j\otimes z_{ij}$ for all $1 \leq i \leq q$. Moreover, since $e,f$ are $L$-traceable, our previous discussion shows that $F(\kk (D^{e,e})^{\pm 1})=(\kk (D^{e,e})^{\pm 1})\cotens_{\mathcal{SL}_m(e)}T=\kk (D^{f,f})^{\pm 1}$. So there are invertible elements $d^{\pm 1}\in T$ such that $(D^{f,f})^{\pm 1}=(D^{e,e})^{\pm 1}\otimes d^{\pm 1}$. It is straightforward to check that the right $\mathcal{GL}_m(f)$-comodule structure on $T$ are 
\begin{gather*}
\rho_r(y_{ji})=\sum_{k=1}^q y_{jk}\otimes a^{f,f}_{ki},\quad \rho_r(z_{ij})=\sum_{k=1}^q z_{kj}\otimes b^{f,f}_{ik},\text{ for } 1\leq j\leq p, 1\leq i\leq q,\quad \rho_r(d^{\pm 1})=d^{\pm 1}\otimes (D^{f,f})^{\pm 1}. 
\end{gather*}

Applying the natural isomorphism of monoidal equivalences $F(-)=(-)\cotens_{\mathcal{GL}_m(e)}T$ to \eqref{Tf}, we get the following commutative diagram 
\[
\xymatrix{
\kk (D^{f,f})^{-1}=F(\kk(D^{e,e})^{-1})\ar[r]^-{F(e)}\ar[d]_-{\cong}\ar@/^25.0pt/[rrr]^-{f} & F(U_e^{\otimes m})\ar[r]^-{\xi^{-1}}\ar[d]_-{\cong}& F(U_e)^{\otimes m}\ar[r]^-{=}\ar[d]_-{\cong } & V_f^{\otimes m}\ar[dl]_-{\cong}\\
 \kk (D^{e,e})^{-1} \cotens_{\mathcal{GL}_m(e)}T\ar[r]^-{e\otimes \id} & (U_e^{\otimes m})\cotens_{\mathcal{GL}_m(e)}T& (U_e\cotens_{\mathcal{GL}_m(e)}T)^{\otimes m}. \ar[l]_-{\xi}& 
}
\]
By a straightforward element tracing as in the proof of \Cref{SQSMTE}, we obtain 
the relations in $T$ 
\[\sum_{1\leq i_1,\ldots,i_m\leq q} f_{i_1\cdots i_m}y_{j_1i_1}\cdots y_{j_mi_m}~=~e_{j_1\cdots j_m}d^{-1},\quad \text{$\forall$ $1\leq j_1,\cdots,j_m\leq p$.}\] 

By doing similarly on $f^*$, we obtain the following commutative diagram 
\[
\xymatrix{
(V_f^*)^{\otimes m}= F(U_e^*)^{\otimes m}\ar@/^34.0pt/[rrrr]^-{f^*} \ar[r]^-{\xi}\ar[d]_-{\cong} & F((U_e^*)^{\otimes m})\ar[r]^-{F(\cong)}\ar[d]_-{\cong}  & F((U_e^{\otimes m})^*)\ar[d]_-{\cong} \ar[r]^-{F(e^*)}\ar[d]_-{\cong}   & F(\kk D^{e,e})\ar[d]_-{\cong} \ar[r]^-{\cong}& \kk D^{f,f}\ar[dl]^-{\cong} \\
 (U_e^*\cotens_{\mathcal{GL}_m(e)}T)^{\otimes m}\ar[r]^-{\xi} & (U_e^*)^{\otimes m} \cotens_{\mathcal{GL}_m(e)}T\ar[r]^-{(\cong) \otimes \id} & (U_e^{\otimes m})^* \cotens_{\mathcal{GL}_m(e)}T\ar[r]^-{e^*\otimes \id}& \kk D^{e,e}\cotens_{\mathcal{GL}_m(e)}T
}
\]
which yields, for all $1 \leq j_1,\cdots, j_m \leq q$,
\begin{align*}
    D^{e,e}\otimes f_{j_1\cdots j_m}d&\,=f_{j_1\cdots j_m}D^{f,f}\\
    &\,=f^*(v^{j_m}\otimes \cdots \otimes v^{j_1})\\
    &\,=f^*\left((\sum_{i_m} u^{i_m}\otimes z_{j_mi_m})\otimes \cdots \otimes(\sum_{i_1} u^{i_1}\otimes z_{j_1i_1})\right)\\
    &\,=\sum_{1\leq i_1,\ldots,i_m\leq p}e^*(u^{i_m}\otimes \cdots \otimes u^{i_1})\otimes z_{j_mi_m}\cdots z_{j_1i_1}\\
    &\, =D^{e,e} \otimes \left(\sum_{1\leq i_1,\ldots,i_m\leq p} e_{i_1\cdots i_m}z_{j_mi_m}\cdots z_{j_1i_1}\right),
\end{align*}
and similarly we obtain
\[
\xymatrix{
\quad V_f^*\otimes V_f\ar[r]^-{=}\ar[dr]_-{\cong}\ar@/^40.0pt/[rrrr]^-{\rm ev}  & F(U_e^*)\otimes F(U_e)\ar[r]^-{\xi}\ar[d]_-{\cong} & F(U_e^*\otimes U_e)\ar[r]^-{F({\rm ev})} \ar[d]_-{\cong}& F(\kk)\ar[r]^-{=}\ar[d]_-{\cong} & \kk\ar[dl]^-{\cong}\\
& (U_e^*\cotens_{\mathcal{GL}_m(e)}T)\otimes (U_e\cotens_{\mathcal{GL}_m(e)}T)\ar[r]^-{\xi}& (U_e^*\otimes U_e)\cotens_{\mathcal{GL}_m(e)}T\ar[r]^-{{\rm ev}\otimes \id} & \kk \cotens_{\mathcal{GL}_m(e)}T
}
\]
which yields, for all $1 \leq i,j \leq q$,
\begin{align*}
    \delta_{ij}=\langle v^i,v_j\rangle= \Big\langle \sum_k u^k\otimes z_{ik}, \sum u_l\otimes y_{lj} \Big\rangle=\sum_{k,l} \langle u^k,u_l\rangle z_{ik}y_{lj}=\sum_{k} z_{ik}y_{kj}.
\end{align*}

In summary, the elements $\{y_{ji}\},\{z_{ij}\},\ d^{\pm 1}$ satisfy the following relations in $T$:
\begin{align*}
\sum_{1\leq i_1,\ldots,i_m\leq q} f_{i_1\cdots i_m}y_{j_1i_1}\cdots y_{j_mi_m}~&=~e_{j_1\cdots j_m}d^{-1},\quad \text{$\forall$ $1\leq j_1,\cdots,j_m\leq p$};\\
  \sum_{1\leq i_1,\ldots,i_m\leq p} e_{i_1\cdots i_m}z_{j_mi_m}\cdots z_{j_1i_1}~&=~f_{j_1\cdots j_m}d,\quad \text{$\forall$ $1\leq j_1,\cdots,j_m\leq q$};  \\
  \sum_{1\leq k\leq p} z_{ik}y_{kj}~&=~\delta_{ij},\quad \text{ $\forall$ $1\leq i,j\leq q$. }
\end{align*}
According to \eqref{eq:alg}, there is a well-defined right $\mathcal{GL}_m(f)$-comodule algebra map $\pi: \mathcal{GL}_m(e,f)\to T$ via $a_{ji}^{e,f}\mapsto y_{ji}$, $b_{ij}^{e,f}\mapsto z_{ij}$ and $(D^{e,f})^{\pm 1}\mapsto d^{\pm 1}$, for all $1 \leq j \leq p, 1 \leq i \leq q$. Moreover, since both $\mathcal{GL}_m(e,f)$ and $T$ are right Galois objects over $\mathcal{GL}_m(f)$, $\pi$ is indeed an isomorphism \cite[Proposition 1.6]{Bichon2014}. Consequently, $\mathcal{GL}_m(e,f)\neq 0$. It is clear that $\pi$ is also a left $\mathcal{GL}_m(e)$-comodule algebra isomorphism since the left $\mathcal{GL}_m(e)$-comodule structure on $T$ is induced by 
\begin{gather*}
\rho_l(y_{ji})=\sum_{k=1}^p a^{e,e}_{jk}\otimes y_{ki},\quad \rho_l(z_{ij})=\sum_{k=1}^p b^{e,e}_{kj}\otimes z_{ik}, \text{ for } 1\leq j \leq p, 1\leq i \leq q,\quad \rho_l(d^{\pm 1})=(D^{e,e})^{\pm 1}\otimes d^{\pm 1}. 
\end{gather*}
The last assertion about $(-)\cotens_{\mathcal{GL}_m(e)}\mathcal{GL}_m(e,f)$ follows directly from the above argument.
\end{proof}

The following result can be proved similarly to that in \cite[Theorem 5.33]{CWW2019}.

\begin{lemma}
\label{SQSAS}
Let $e\in U^{\otimes m},f\in V^{\otimes m}$ be $L$-traceable nondegenerate twisted superpotentials with integers $2 \leq L\leq N\leq m$, then we have isomorphisms of universal quantum groups
    \[
    \mathcal{GL}_m(e)~\cong~\underline{\rm aut}^r(A(e,N)) \qquad \text{ and } \qquad \mathcal{GL}_m(f)~\cong~\underline{\rm aut}^r(A(f,N)),
    \]
    where the universal coactions of $\mathcal{GL}_m(e)$ and $\mathcal{GL}_m(f)$ on $A(e,N)$ and $A(f,N)$ are induced by those on $U_e$ and $V_f$, respectively. In particular, the following are equivalent.
    \begin{enumerate}
    \item[(1)] $\mathcal{GL}_m(e,f)\neq 0$.
        \item[(2)] $A(e,N)$ and $A(f,N)$ are $\GLQS$ equivalent.
        \item[(3)] $A(e,N)$ and $A(f,N)$ are QS equivalent. 
    \end{enumerate} 
\end{lemma}

\begin{proof}
For the first part of the lemma, it suffices to show that the Hopf coaction of $\mathcal{GL}_m(e)$ on $A(e,N)$ has the universal property: for any Hopf coaction of $H$ on $A(e,N)$ while preserving its grading, there is a unique Hopf algebra map $\varphi: \mathcal{GL}_m(e)\to H$ such that the $H$-coaction on $A(e,N)$ factors through that of $\mathcal{GL}_m(e)$ via $\varphi$. Since $H$ coacts on $A(e,N)=TU/\langle\partial^{m-N}(\kk e)\rangle$ while preserving its grading, $\partial^{m-N}(\kk e)\subset U^{\otimes N}$ is a subcomodule over $H$. Consequently, $\partial^{m-L}(\kk e)=\partial^{N-L}(\partial^{m-N}(\kk e))\subset U^{\otimes L}$ is again an $H$-subcomodule. Because $e\in U^{\otimes m}$ is $L$-traceable, the natural inclusion 
\[
\kk e=\bigcap_{i+j=m-L} U^{\otimes i}\otimes \partial^{m-L}(\kk e)\otimes U^{\otimes j}\subset U^{\otimes m}
\] 
is $H$-collinear. Then we can apply the universal property of $\mathcal{GL}_m(e)$ regarding Hopf coactions on the above embedding $\kk e\hookrightarrow U^{\otimes m}$ (see \cite[Proposition 5.28]{CWW2019}) to get a unique Hopf algebra map $\varphi: \mathcal{GL}_m(e)\to H$ such that the $H$-coaction on $U$ factors through that of $\mathcal{GL}_m(e)$. Consequently, $\varphi$ is the unique Hopf algebra map such that the $H$-coaction on $A(e,N)$ factors through that of $\mathcal{GL}_m(e)$. Thus, by the universal property of Manin's universal quantum groups, we get the isomorphisms between these universal quantum groups $\mathcal{GL}_m(e)~\cong~\underline{\rm aut}^r(A(e,N))$. The same argument holds for the assertion involving $f$.

For the second part of the lemma, (1) $\Rightarrow$ (2) comes from  \Cref{SQSMTE}; (2) $\Leftrightarrow$ (3) is due to the first part of the lemma, and (2) $\Rightarrow$ (1) follows from \Cref{SQSBG}.
\end{proof}

When $m=2$, every nondegenerate twisted superpotential $f\in V^{\otimes 2}$ is given by some invertible matrix $\mathbb F=(f_{ij})_{1\leq i,j\leq q}\in \GL(V)$ such that $f=\sum_{1\leq i,j\leq q} f_{ij}v_i\otimes v_j$ where $\{v_1,\ldots,v_q\}$ is a basis for $V$. We know $A(f,2)\cong \kk\langle v_1,v_2,\cdots,v_q\rangle/(\sum_{1\leq i,j\leq q} f_{ij}v_iv_j)$ is an AS-regular algebra of global dimension two by \cite{Zhang1996}. As a consequence, every nondegenerate twisted superpotential in $V^{\otimes 2}$ is 2-traceable, and the cogroupoid $\mathcal{GL}_2$ is connected \cite[Theorem B]{HNUVVW2024-3}. For any integer $m \geq 2$, the following theorem shows that $\mathcal{GL}_m$ has a connected component consisting of all Koszul AS-regular superpotentials. Recall the notation ${\rm AS}(m,2,q)$ given in \eqref{ASclass} and note that by \cite[Proposition 2.12]{MS2016}, any superpotential $f$ which gives rise to a Koszul AS-regular superpotential algebra $A(f,2)$ is $2$-traceable.

\begin{thm}
\label{ASSuper}
    Let $m\ge 2$ be any integer and $f\in V^{\otimes m}$ a nondegenerate twisted superpotential with $\dim V = q$. The following are equivalent.
    \begin{enumerate}
        \item[(1)] $\kk f\in {\rm AS}(m,2,q)$;
        \item[(2)] $\mathcal{GL}_m(e,f)\neq 0$ for all $\kk e\in\bigsqcup_{p\geq 2} {\rm AS}(m,2,p)$; 
        \item[(3)] $\mathcal{GL}_m(e,f)\neq 0$ for some $\kk e\in\bigsqcup_{p\geq 2} {\rm AS}(m,2,p)$. 
    \end{enumerate}
    Consequently, all superpotentials in $\bigsqcup_{p\geq 2} {\rm AS}(m,2,p)$ giving Koszul AS-regular algebras of dimension $m$ form a connected component of the cogroupoid $\mathcal{GL}_m$. 
\end{thm}

\begin{proof}
 (1) $\Rightarrow$ (2): Suppose $f$ is Koszul AS-regular. Let $e\in U^{\otimes m}$ be any Koszul AS-regular twisted superpotential with $\dim U=p$. By \Cref{NKAS}, we know $A(e,2)$ and $A(f,2)$ are Koszul AS-regular algebras of the same global dimension $m$. First, assume $\kk$ to be algebraically closed. Then, by \cite[Theorem 5.1.1]{HNUVVW2024-2}, $A(e,2)$ and $A(f,2)$ are QS equivalent since they share the same global dimension. Since $e,f$ are $2$-traceable, by \Cref{SQSAS}, we have $\mathcal{GL}_m(e,f)\neq 0$. If $\kk$ is not algebraically closed, consider $A(f,2) \otimes \overline{\kk}$, which is now AS-regular over an algebraically closed field. By the above argument, $\mc{GL}_m(e,f) \otimes \overline{\kk}$ is nonzero, which then implies that $\mc{GL}_m(e,f) \not \neq 0$.

 (2) $\Rightarrow$ (3): is routine. 
 
 (3) $\Rightarrow$ (1): By the assumption of (3) and \Cref{SQSMTE}, $A(e,2)$ and $A(f,2)$ are $\GLQS$ equivalent for some Koszul AS-regular superpotential $e \in U^{\otimes m}$. Let $s\in W^{\otimes m}$ be the twisted superpotential such that $A(s,2)\cong \kk[x_1,\ldots,x_m]$. By the discussion in (1) $\Rightarrow$ (2), we know $\mathcal{GL}_m(s,e)\neq 0$. So $A(s,2)$ and $A(e,2)$ are $\GLQS$ equivalent by \Cref{SQSMTE}. This implies that $A(s,2)$ and $A(f,2)$ are $\GLQS$ equivalent. Note that $A(s,2)$ is just a noetherian polynomial algebra and $\mathcal{GL}_m(s)$ has bijective antipode \cite[Corollary 5.27]{CWW2019}. Then we can conclude that $A(f,2)$ is Koszul AS-regular by \Cref{PQS}. The last statement is clear.
\end{proof}

In the following result, we demonstrate that the class of Koszul AS-regular algebras of the same dimension $m$ forms a single quantum-symmetric class that includes the polynomial algebra in $m$ variables. This observation provides an alternative characterization of Koszul AS-regular algebras as quantum deformations of projective spaces from a quantum-symmetric perspective. This result extends \cite[Theorem H]{HNUVVW2024-2}, where we relax the assumption regarding the base field, which was previously assumed to be algebraically closed due to the categorical approach of the Tannaka--Krein formalism in \cite{RVdB2017}.

\begin{thm}
\label{AS2cocyle}
Let $A$ be a Koszul AS-regular algebra of finite global dimension and $B$ any connected graded algebra finitely generated in degree one. Then $A$ and $B$ are QS equivalent if and only if $B$ is a Koszul AS-regular algebra of the same global dimension as $A$. In this case, $A$ and $B$ are 2-cocycle twists of each other if and only if they have the same Hilbert series.   
\end{thm}

\begin{proof}
Let $A$ be Koszul AS-regular of finite global dimension $m$. By \Cref{NKAS}, we may write $A=A(e,2)$ for some $\kk e\in {\rm AS}(m,2,p)$ . 

$\Rightarrow$: Supppose that $A,B$ are QS equivalent. Let $\kk s\in {\rm AS}(m,2,m)$ such that $A(s,2)=\kk[x_1,\ldots,x_m]$. By \Cref{ASSuper}, $\mathcal{GL}_m(s,e)\neq 0$. Hence, \Cref{SQSAS} implies that $A(s,2)$ and $A(e,2)$ are QS equivalent and,  so are $A(s,2)$ and $B$. Note that $A(s,2)$ is the noetherian polynomial algebra and $\mathcal{GL}_m(s)$ has a bijective antipode \cite[Corollary 5.27]{CWW2019}. Thus, we can conclude that $B$ is Koszul AS-regular of global dimension $m$, by \Cref{PQS}. 

$\Leftarrow$: Conversely, suppose that $B$ is Koszul AS-regular of global dimension $m$. Again by \Cref{NKAS}, we can write $B=A(f,2)$ with $\mathcal{GL}_m(e,f)\neq 0$ by \Cref{ASSuper}. Hence, $A$ and $B$ are QS equivalent, by \Cref{SQSAS}.

For the last statement of the theorem, assume that $A$ and $B$ are two Koszul AS-regular algebras of the same global dimension $m$. Without loss of generality, we can write $A=A(e,2)$ and $B=A(f,2)$ for some $\kk e\in {\rm AS}(m,2,p)$ and $\kk f\in {\rm AS}(m,2,q)$. By \Cref{ASSuper}, $\mathcal{GL}_m(e,f)\neq 0$ and $F(-)=(-)\cotens_{\mathcal{GL}_m(e)}\mathcal{GL}_m(e,f)$ induces a $\GLQS$ equivalence between $A$ and $B$. By writing $A=TU_e/(R_e)$, where $R_e=\partial^{m-2}(\kk e)$ is the vector space spanned by the relation which define $A$ as in \Cref{lem:comoduleGL}, one sees that there are comodules $M_A(r_i)$ over $\mathcal{GL}_m(e)$ where $M_A(r_1)=U_e$ and  $M_A(r_k)=\bigcap_{i+j+2=k} U_e^{\otimes i}\otimes R_e\otimes U_e^{\otimes j}$ for $2\leq k\leq  m$, where $M_A(r_m)$ is the invertible comodule given by the grouplike element $(D^{e,e})^{-1}$. We use the same notation $M_B(r_i)$ with $1 \leq i\leq m$ for the corresponding comodules over $\mathcal{GL}_m(f)$. Since $F$ sends $A$ to $B$ as comodule algebras, it follows that $F(M_A(r_i))=M_B(r_i)$, for all $1 \leq i\leq m$. We let $A'$ be a possible base field extension of $A$ and write $A'=A\otimes_\kk \overline{\kk}=A(e',2)$ and $B'=B\otimes_\kk \overline{\kk}=A(f',2)$, where $e'$ and $f'$ are the corresponding nondegenerate twisted superpotentials for $e,f$ after the base field extension. In particular, $M_{A'}(r_i)=M_A(r_i)\otimes_\kk \overline{\kk}$ and $M_{B'}(r_i)=M_B(r_i)\otimes_\kk \overline{\kk}$ become comodules over the corresponding universal quantum groups $\mathcal{GL}_m(e')=\mathcal{GL}_m(e)\otimes_\kk \overline{\kk}$ and $\mathcal{GL}_m(f')=\mathcal{GL}_m(f)\otimes_\kk \overline{\kk}$. Moreover, the $\mathcal{GL}_m(e')$-$\mathcal{GL}_m(f')$-bi-Galois object 
$\mathcal{GL}_m(e',f')=\mathcal{GL}_m(e,f)\otimes_\kk\overline{\kk}$ implies a  $\GLQS$ equivalence, and hence a QS equivalence by \Cref{SQSAS}, between $A'$ and $B'$.

If $A,B$ have the same Hilbert series, then we have 
\begin{equation}\label{eq:dimeq}
  \dim_{\overline{\kk}} M_{A'}(r_i)= \dim_\kk M_A(r_i)=\dim_\kk (A^!)_i=\dim_\kk (B^!)_i=\dim_\kk M_B(r_i) =\dim_{\overline{\kk}} M_{B'}(r_i)
\end{equation} 
for all $1\leq i\leq m$. Since $\mathcal{GL}_{m}(e')$ is the Manin's universal quantum group for $A'=A(e',2)$, by \cite[Corollary 15.43]{Manin2018}, we know that the Grothendieck ring of $\mathcal{GL}_m(e')$  is the free $\mathbb Z$-ring 
\[
\mathbb Z\langle M_{A'}(r_1),\ldots, M_{A'}(r_{m-1}), M_{A'}(r_m)^{\pm 1}\rangle,
\]
and the same is true for $\mathcal{GL}_m(f')$. Consequently, the monoidal equivalence 
\[F': {\rm comod}\left(\mathcal{GL}_m(e')\right)\to\comod\left(\mathcal{GL}_m(f')\right)\] induces a ring isomorphism $F'$ (by abuse of notation) between the Grothendieck rings of $\mathcal{GL}_m(e')$ and $\mathcal{GL}_m(f')$ by sending $M_{A'}(r_i)$ to $M_{B'}(r_i)$ for all $1\leq i\leq m$. For any finite-dimensional comodule $W$ over the original $\mathcal{GL}_m(e)$, we can consider $W'=W\otimes_\kk\overline{\kk}$ as an element in the Grothendieck ring of $\mathcal{GL}_m(e')$ and write it as a $\mathbb Z$-linear combination of some tensor products of the free generators  $M_{A'}(r_1),\ldots,M_{A'}(r_m)^{\pm 1}$. By applying $F'$, one can similarly write $F'(W')$ as the $\mathbb Z$-linear combination of the tensor products now in terms of $M_{B'}(r_1),\ldots,M_{B'}(r_m)^{\pm 1}$ in the Grothendieck ring of $\mathcal{GL}_m(f')$.  This implies that $\dim_\kk F(W)=\dim_{\overline{\kk}} F'(W')=\dim_{\overline{\kk}}W'=\dim_\kk W$ due to \eqref{eq:dimeq} and hence $F$ preserves dimensions. As a consequence, $\mathcal{GL}_m(f)$ is a 2-cocycle twist of $\mathcal{GL}_m(e)$ (see \cite[Theorems 1.8 and 1.17]{Bichon2014} and \cite[Proposition 4.2.2]{PavelGelaki2001}), where $F$ is naturally isomorphic to the monoidal equivalence given by the 2-cocycle twist. Since $F$ maps the comodule algebra $A$ over $\mathcal{GL}_m(e)$ to the comodule algebra $B$ over $\mathcal{GL}_m(f)$, we know $B$ is a 2-cocycle twist of $A$. Conversely, if $A$ and $B$ are 2-cocycle twists of each other, they share the same Hilbert series by construction.
\end{proof}

\begin{remark}
Recall that we do not assume AS-regular algebras have finite GK-dimension in this paper; the AS-regular algebras with infinite GK-dimension will never be 2-cocycle twists of AS-regular algebras with finite GK-dimension. On the other hand, the Hilbert series of the polynomial ring in $m$-variables is $\frac{1}{(1-t)^m}$, and it is an open question whether or not all Koszul AS-regular algebras of global dimension $m$ of finite GK-dimension share this Hilbert series (see \cite[Question 3.3(5)]{Rogalski2024}). Our result in \Cref{AS2cocyle} shows that Koszul AS-regular algebras have Hilbert series $\frac{1}{(1-t)^m}$ if and only if they are in the same QS equivalence class with the polynomial ring in $m$-variables and this equivalence is given by a 2-cocycle twist.
\end{remark}

\begin{lemma}
\label{isobi-Galois}
Let $m \geq 2$ be any integer. Let $e,e'\in U^{\otimes m}$ and $f,f'\in V^{\otimes m}$ be nondegenerate twisted superpotentials. Suppose there are $\phi \in \GL(U)$ and $\psi \in \GL(V)$ such that $\phi^{\otimes m}(\kk e)=\kk e'$ and $\psi^{\otimes m}(\kk f)=\kk f'$. Then $\mathcal{GL}_m(e,f)$ and $\mathcal{GL}_m(e',f')$ are isomorphic as $\kk$-algebras. In particular, $\mathcal{GL}_m(e,f)\cong \mathcal{GL}_m(e,f')$ as left $\mathcal{GL}_m(e)$-comodule algebras and $\mathcal{GL}_m(e,f)\cong \mathcal{GL}_m(e',f)$ as right $\mathcal{GL}_m(f)$-comodule algebras. 
\end{lemma}
\begin{proof}
Choose bases $\{u_1,\ldots,u_p\}$ for $U$ and $\{v_1,\ldots,v_q\}$ for $V$. Write $\phi(u_i)=\sum_{1\leq j\leq p} u_j\phi_{ji}$ and $\psi(v_i)=\sum_{1\leq j\leq q} v_j\psi_{ji}$ for invertible matrices $\phi=(\phi_{ij})\in \GL(U)$ and $\psi=(\psi_{ij})\in \GL(V)$. By our assumption, we also have $\phi^{\otimes m}(e)=\alpha e'$ and $\psi^{\otimes m}(f)=\beta f'$ for some nonzero scalars $\alpha,\beta$. According to \eqref{eq:alg}, we can check that the following map is an algebra isomorphism from $\mathcal{GL}_m(e,f)$ to $\mathcal{GL}_m(e',f')$:
\[
\mathbb A^{e,f}\mapsto \phi \mathbb A^{e',f'} \psi,\qquad \mathbb B^{e,f}\mapsto \psi^{-1}\mathbb B^{e',f'}\phi^{-1},\qquad (D^{e,f})^{\pm 1}\mapsto (D^{e',f'}/\alpha\beta)^{\pm 1}.
\]
For the last in particular statement, whenever $e=e'$ or $f=f'$, we can choose either $\phi=\mathbb I_{p\times p}$ or $\psi=\mathbb I_{q\times q}$ and the above isomorphism becomes a left $\mathcal{GL}_m(e)$-comodule map or a right $\mathcal{GL}_m(f)$-comodule map.  
\end{proof}

We now reformulate the classification of Koszul AS-regular algebras in terms of the classification of fiber functors from the comodule categories over Manin's universal quantum groups of the polynomial algebras.

\begin{thm}
\label{IsoAS}
For any integer $m \geq 2$, the following are in one-to-one correspondence with each other.
\begin{enumerate}
    \item[(1)] Isomorphism classes of Koszul AS-regular algebras of global dimension $m$.
    \item[(2)] The set of left $\underline{\rm aut}^r(\kk[x_1,\ldots,x_m])$-Galois objects up to isomorphism.
    \item[(3)] The set of fiber functors from $\comod_{\fd}(\underline{\rm aut}^r(\kk[x_1,\ldots,x_m]))$ to ${\rm Vec}(\kk)$ up to natural isomorphism.     
    \item[(4)] Disjoint union of orbits $\bigsqcup_{q\ge 2} {\rm AS}(m,2,q)/\PGL_q(\kk)$, where ${\rm AS}(m,2,q)$ is defined in \eqref{ASclass}.
\end{enumerate}
\end{thm}

\begin{proof}
    (1) $\Leftrightarrow$ (4): Since Koszul AS-regular algebras are connected graded algebras that are finitely generated in degree one, we know that two of them are isomorphic as associative algebras if and only if they are isomorphic as graded algebras \cite[Theorem 1]{BZ2017}. Then by \Cref{NKAS}, it suffices to show two superpotential algebras $A(e,2)$ and $A(f,2)$ given by two nondegenerate twisted superpotentials $\kk e\in {\rm AS}(m,N,p)$ and $\kk f\in {\rm AS}(m,N,q)$ are isomorphic as graded algebras if and only if there is a linear isomorphism $\phi: U\cong V$ such that $\phi^{\otimes m}(\kk e)=\kk f$, which is well-known \cite[Theorem 2.11]{MU2019}.

    (2) $\Leftrightarrow$ (3): This follows from Ulbrich's theorem \cite{Ulbrich1987}, see also \cite[Theorem 1.14]{Bichon2014}.

(1),(4) $\Leftrightarrow$ (2),(3): Let $\kk s\in {\rm AS}(m,2,m)$ such that $A(s,2) \cong \kk[x_1,\ldots,x_m]$. By \Cref{SQSAS}, we have $\underline{\rm aut}^r(\kk[x_1,\ldots,x_m])\cong \mathcal{GL}_m(s)$. For every nondegenerate twisted superpotential $\kk f\in {\rm AS}(m,2,q)$, we know that $\mathcal{GL}_m(s,f)\neq 0$ yields a $\mathcal{GL}_m(s)$-$\mathcal{GL}_m(f)$-bi-Galois object by \Cref{ASSuper}. Hence, we obtain a fiber functor from $\comod_{\fd}(\underline{\rm aut}^r(\kk[x_1,\ldots,x_m]))$ to ${\rm Vec}(\kk)$ by composing $(-)\cotens_{\mathcal{GL}_m(s)}\mathcal{GL}_m(s,f)$ with the forgetful functor $\comod_{\fd}(\mathcal{GL}_m(f)) \to {\rm Vec}(\kk)$. Conversely, suppose $T$ is a left $\mathcal{GL}_m(s)$-Galois object. By \Cref{Schauenburg}, there is a unique Hopf algebra $K$ such that $T$ becomes a $\mathcal{GL}_m(s)$-$K$-bi-Galois object and  $(-)\cotens_{\mathcal{GL}_m(s)}T$ induces a monoidal equivalence between $\comod(\mathcal{GL}_m(s))$ and $\comod(K)$. Denote by $B=A(s,2)\cotens_{\mathcal{GL}_m(s)}T$ the comodule algebra in $\comod(K)$, which is Koszul AS-regular of global dimension $m$ by \cite[Corollary 3.2.7 \& Theorem 3.2.10]{HNUVVW2024-2}. We write $B=A(f,2)$ for some nondegenerate twisted superpotential $\kk f\in {\rm AS}(m,2,q)$. Similar to the second part of the proof of \Cref{RQS}, we can show that $K\cong \underline{\rm aut}^r(A(f,2))\cong \mathcal{GL}_m(f)$. Therefore, we know $(-)\cotens_{\mathcal{GL}_m(s)}T$ induces an $\GLQS$ equivalence between $A(s,2)$ and $A(f,2)$. By \Cref{SQSBG}, we have $T\cong \mathcal{GL}_m(s,f)$ as $\mathcal{GL}_m(s)$-$\mathcal{GL}_m(f)$-bicomodule algebras. It remains to show that for two $f\in {\rm AS}(m,2,q)$ and $f'\in {\rm AS}(m,2,q')$, we have $\mathcal{GL}_m(s,f)$ and $\mathcal{GL}_m(s,f')$ are isomorphic as left $\mathcal{GL}_m(s)$-comodule algebras if and only if $A(f,2)$ and $A(f',2)$ are isomorphic as (graded) algebras. Denote by $F(-)=(-)\cotens_{\mathcal{GL}_m(s)}\mathcal{GL}_m(s,f)$ and $F'(-)=(-)\cotens_{\mathcal{GL}_m(s)}\mathcal{GL}_m(s,f')$ the corresponding fiber functors. One direction is clear since $F(A(s,2))=A(f,2)$ and $F'(A(s,2))=A(f',2)$ by \Cref{SQSMTE}. If $F\cong F'$ are isomorphic as monoidal functors, we must have $A(f,2)\cong A(f',2)$ as graded algebras. On the other hand, suppose $A(f,2)=TV/\langle\partial^{m-2}(\kk f)\rangle \cong A(f',2)=TV'/\langle\partial^{m-2}(\kk f')\rangle$ as graded algebras, we have an isomorphism $\phi: V\cong V'$ such that $\phi^{\otimes m}(\kk f)=\kk f'$. Then, $\mathcal{GL}_m(s,f)\cong \mathcal{GL}_m(s,f')$ as left $\mathcal{GL}_m(s)$-comodule algebras by \Cref{isobi-Galois}.
\end{proof}

The following result identifies the isomorphism classes of all Koszul AS-regular algebras that are deformations of polynomial algebras with left cleft Galois objects over Manin's universal quantum groups of the polynomial algebras. 

\begin{Cor}
For any integer $m \geq 2$, the following are in one-to-one correspondence with each other.
\begin{itemize}
    \item[(1)] Isomorphism classes of Koszul AS-regular algebras with Hilbert series $\frac{1}{(1-t)^m}$.
        \item[(2)] The set of left cleft $\underline{\rm aut}^r(\kk[x_1,\ldots,x_m])$-Galois objects up to isomorphism.
    \item[(3)] The set of fiber functors from $\comod_{\fd}(\underline{\rm aut}^r(\kk[x_1,\ldots,x_m]))$ to ${\rm Vec}(\kk)$ preserving vector space dimensions up to natural isomorphism. 
\end{itemize}   
\end{Cor}

\begin{proof}
(2) $\Leftrightarrow$ (3) This is due to Etingof-Gelaki \cite{EG2001}, see also \cite[Theorem 1.17]{Bichon2014}.

(1) $\Leftrightarrow$ (2),(3) By the work of Doi-Takeuchi \cite{DT86} and Blattner-Montgomery \cite{BM1989}, we know that every cleft Galois object is given by some 2-cocycle, whose corresponding fiber functor sends comodule algebras to their 2-cocycle twists. Then the sub-correspondence follows from \Cref{AS2cocyle}.
\end{proof}

\begin{remark}
It is known that any AS-regular algebra with Hilbert series $\frac{1}{(1-t)^m}$ is Koszul \cite[Theorem 2.2]{BC2001} and then can be written as $A(f,2)$ for some $\kk f\in {\rm AS}(m,2,m)$ (cf. \Cref{NKAS}).
\end{remark}

\subsection{Some invariants under $\SLQS$ equivalence}
\label{subsec:SLQS}

We can obtain similar results for $\mathcal{SL}_m$ regarding $\SLQS$ equivalence by slightly modifying the arguments for $\mathcal{GL}_m$.

\begin{thm}
\label{thm:SLQS}
Let $2\leq N\leq m$ be integers and $A(e,N)$ and $A(f,N)$ be two superpotential algebras for some nondegenerate twisted superpotentials $e\in U^{\otimes m}$ and $f\in V^{\otimes m}$, respectively. 
\begin{itemize}
    \item[(1)]If $\mathcal{SL}_m(e,f)\neq 0$, then $A(e,N)$ and $A(f,N)$ are $\SLQS$ equivalent. In particular, they are $\GLQS$ equivalent. 
    \item[(2)]If $A(e,N)$ and $A(f,N)$ are $\SLQS$ equivalent for some $L$-traceable $e,f$ with $2\leq L\leq N\leq m$, then $\mathcal{SL}_m(e,f)\neq 0$. 
\end{itemize} 
\end{thm}

\begin{proof}
    (1): Since $\mathcal{GL}_m(e,f)=\mathcal{SL}_m(e,f)/(D^{e,f}-1)$, $\mathcal{GL}_m(e,f)\neq 0$ provided $\mathcal{SL}_m(e,f)\neq 0$. Then (1) follows from an analogous proof of \Cref{SQSMTE} for both $\SLQS$ equivalence and $\GLQS$ equivalence conclusions. 
    
    (2): Suppose $A(e,N)$ and $A(f,N)$ are $\SLQS$ equivalent for some $L$-traceable $e,f$. We can apply a similar argument in \Cref{SQSBG} to conclude that  $\mathcal{SL}_m(e,f)\neq 0$.
\end{proof}

Recall that by \Cref{L:copi}, the cogroupoid $\mathcal{SL}_m$ is copivotal such that 
for any nondegnerate twisted superpotential $e\in U^{\otimes m}$ with associated invertible matrix $\mathbb P$, the copivotal character $\Phi_e: \mathcal{SL}_m(e)\to \kk$ is given by $\Phi_e(\mathbb A^{e,e})=\mathbb P^T, \Phi_e(\mathbb B^{e,e})=\mathbb P^{-T}$. Moreover, for another nondegenerate twisted superpotential $f\in V^{\otimes m}$, we have $S_{f,e}\circ S_{e,f}=\Phi_e^{-1}* \id *\Phi_f$ with $\Phi_e^{-1}=\Phi_e\circ S_{e,e}$ denotes the convolution inverse of $\Phi_e$. Let $W\in \comod_{\fd}(\mathcal{SL}_m(e))$. Choose a basis $\{w_1,\ldots,w_r\}$ for $W$ and write the $\mathcal{SL}_m(e)$-comodule structure map on $W$ as $\rho_W(w_i)=\sum_j w_j\otimes h_{ji}$ for some $h_{ji}\in \mathcal{SL}_m(e)$. 

\begin{lemma}[{\cite[Proposition 2.9]{Bichon2001}}]
Retain the above notation. For any $W\in \comod_{\fd}(\mathcal{SL}_m(e))$ of dimension $r$, there is a natural isomorphism $j_W: W\to W^{**}$ given by 
\begin{align}
\label{jW}
j_W(w_i)~=~(\id \otimes \Phi_e^{-1})(\rho(w_i))~=~\sum_j w_j\,\Phi_e^{-1}(h_{ji}), \quad \text{for } 1 \leq i,j \leq r.
\end{align}
\end{lemma}

Moreover, one can extend $j: \id \to (-)^{**}$ to a monoidal isomorphism and such $j$ is said to be a \emph{pivotal structure} on $\comod_{\fd}(\mathcal{SL}_m(e))$. Using this pivotal structure on $\comod_{\fd}(\mathcal{SL}_m(e))$, we recollect the following notions (see e.g., \cite[Section 4.7]{EGNO2015}), which are shown to be invariant given the non-vanishing of $\mathcal{SL}_m(e,f)$, see \Cref{prop:SameH} and \Cref{equivm2}.

\begin{defn}
\label{defn:qdim}
Let $2\leq N\leq m$ be integers. Suppose $e\in U^{\otimes m}$ is a nondegenerate twisted superpotential. Let $W\in \comod_{\fd}(\mathcal{SL}_m(e))$.
\begin{itemize}
    \item[(1)] The \emph{quantum dimension} of $W$ is defined to be 
    \[d_e(W) \coloneqq \left (\kk\xrightarrow{\rm coev} W\otimes W^*\xrightarrow{j_W\otimes \id} W^{**}\otimes W^*\xrightarrow{\rm ev} \kk\right).\]
    \item[(2)] The quantum group $\mathcal{SL}_m(e)$ is  \emph{cospherical} if $d_e(W)=d_e(W^*)$ for any $W\in \comod_{\fd}(\mathcal{SL}_m(e))$.
    \item[(3)] The \emph{quantum Hilbert series} of the superpotential algebra $A(e,N)=\bigoplus_{i\ge 0}A_i$ is defined to be 
\[q_{A(e,N)}(t) \coloneqq \sum_{i\ge 0} d_e(A_i)\,t^i.\]
\end{itemize}
\end{defn}

For the superpotential algebra $A(e,N)=TU/\langle\partial^{m-N}(\kk e)\rangle$, we define 
\[
W_i \coloneqq \begin{cases}
    U^{\otimes i}, &  0\leq i\leq N-1\\
    \bigcap_{s+t=i-N} U^{\otimes s}\otimes \partial^{m-N}(\kk e) \otimes U^{\otimes t}\subseteq U^{\otimes i}, & i\ge N, 
\end{cases}
\]
for any $i\ge 0$ and define   
\[
\rho(i) \coloneqq 
\begin{cases}
    jN, & i=2j\\
    jN+1,  & i=2j+1.
\end{cases}
\]
Suppose $A(e,N)$ is $N$-Koszul AS-regular. By \cite[Theorem 2.4]{BM2006} (also see \cite[Theorem 2.8]{MS2016}), the minimal graded free resolution of the trivial $A(e,N)$-module $\kk$ is given by:
\begin{align}
\label{resolution}
     0\to A(e,N)\otimes W_{\rho(d)}[-\rho(d)]\to \cdots \to A(e,N)\otimes W_{\rho(1)}[-\rho(1)]\to A(e,N)\to \kk,
\end{align}
where $d=2(m-1)/N+1$ is the global dimension of $A(e,N)$. 

\begin{lemma}\label{Le:qd}
  Retain the above notation. Suppose $A(e,N)$ is $N$-Koszul AS-regular of dimension $d$.
  \begin{itemize}
      \item[(1)] The quantum Hilbert series of $A(e,N)$ is 
  \[
  q_{A(e,N)}(t)~=~\frac{1}{\sum\limits_{0\leq i\leq d} (-1)^i\,d_e(W_{\rho(i)})\, t^{\rho(i)}}.
  \]
  \end{itemize}
Further, assume that $A(e,N)$ is Calabi--Yau of dimension $d$. Then we have:
  \begin{itemize}
      \item[(2)] $d_e(W)=(-1)^{d+1}\dim W$, for any $W\in \comod_{\fd}(\mathcal{SL}_m(e))$ and $\mathcal{SL}_m(e)$ is cospherical; and
\item[(3)] $q_{A(e,N)}(t)=(-1)^{d+1}h_{A(e,N)}(t)$, where $h_{A(e,N)}(t)$ is the Hilbert series of $A(e,N)$. 
  \end{itemize} 
\end{lemma}
\begin{proof}
(1): Since the quantum dimension $d_e$ is additive on exact sequences \cite[Proposition 4.7.5]{EGNO2015}, we can derive the quantum Hilbert series of $A(e,N)$ similar to that of the usual Hilbert series from \eqref{resolution}.

(2) \& (3): When $A(e,N)$ is CY of dimension $d$, \Cref{NKAS} says that $\mathbb P=(-1)^{d+1}\mathbb I$, where $\mathbb P$ is the invertible matrix associated to $e$. Then \Cref{L:copi} implies that the character $\Phi_e$ is given by $\Phi_e(\mathbb A)=(-1)^{d+1}\mathbb I$ and $\Phi_e(\mathbb B)=(-1)^{d+1}\mathbb I$. Thus $\Phi_e^{-1}=\Phi_e=(-1)^{d+1}\varepsilon$. So by \eqref{jW}, we get 
  \begin{align*}
    d_e(W)&\,=\sum_{i} \langle j_W(w_i),w^i\rangle= \sum_{i} \langle \sum_jw_j\Phi_e^{-1}(h_{ji}),w^i\rangle=(-1)^{d+1}\sum_i \langle w_j\varepsilon(h_{ji}),w^i\rangle\\
    &\,=(-1)^{d+1}\sum_i\langle w_i,w^i\rangle=(-1)^{d+1}\dim W.
  \end{align*}
Thus, (2) holds, and (2) clearly implies (3). 
\end{proof}

\begin{proposition}
\label{prop:SameH}
Let $2 \leq N\leq m$ be integers and $e\in U^{\otimes m}$ and $f\in V^{\otimes m}$ be two nondegenerate twisted superpotentials. If $\mathcal{SL}_m(e,f)\neq 0$, then $q_{A(e,N)}(t)=q_{A(f,N)}(t)$.
\end{proposition}

\begin{proof}
By \Cref{thm:SLQS}, there is a monoidal equivalence
\begin{align*}
(F,\xi): \comod(\mathcal{SL}_m(e))\to \comod(\mathcal{SL}_m(f)),
\end{align*} 
sending $A(e,N)$ to $A(f,N)$ as comodule algebras. We show that $F$ preserves the copivotal structures (cf. \cite[\S 1]{Ng-Schauenburg}), that is, the diagram
\begin{align}\label{pivS}
\xymatrix{
F(W)\ar[d]_-{j_{F(W)}}\ar[rr]^-{F(j_W)} && F(W^{**})\ar[d]^-{\widetilde{\xi}}\\
F(W)^{**}\ar[rr]^-{\widetilde{\xi}^*} && F(W^*)^*
}
\end{align}
commutes for any $W\in \comod_{\fd}(\mathcal{SL}_m(e))$, where the natural isomorphism $j_W: W\to W^{**}$ is given in \eqref{jW}, and $\widetilde{\xi}$ is defined in \eqref{xitilde}. We write $F(-)=(-)\cotens_{\mathcal{SL}(e)}\mathcal{SL}_m(e,f)$. Fix a basis $\{w_1,\ldots,w_r\}$ for $W$ and write $\rho(w_i)=\sum_j w_j\otimes h_{ji}^{e,e}$ for $(h_{ji}^{e,e})_{1\leq j,i\leq r}\in \mathcal{SL}_m(e)$ satisfying $\Delta_{e,e}^e(h_{ji}^{e,e})=\sum_k h_{jk}^{e,e}\otimes h_{ki}^{e,e}$. Denote $M=F(W)\in \comod_{\fd}(\mathcal{SL}_m(f))$. Similarly, fix a basis $\{m_1,\ldots,m_t\}$ for $M$ and write $\rho(m_i)=\sum_j m_j\otimes h_{ji}^{f,f}$ for $(h_{ji}^{f,f})_{1\leq j,i\leq t}\in \mathcal{SL}_m(f)$ satisfying $\Delta_{f,f}^f(h_{ji}^{f,f})=\sum_k h_{jk}^{f,f}\otimes h_{ki}^{f,f}$. Thus, we can write $m_i=\sum_j w_j\otimes h_{ji}^{e,f}\in W\cotens_{\mathcal{SL}_m(e)}\mathcal{SL}_m(e,f)$ for some $(h_{ji}^{e,f})_{1\leq j\leq r, 1\leq i\leq t}\in \mathcal{SL}_m(e,f)$ satisfying $\Delta_{e,f}^e(h_{ji}^{e,f})=\sum_k h_{jk}^{e,e}\otimes h_{ki}^{e,f}$ and $\Delta_{e,f}^f(h_{ji}^{e,f})=\sum_k h_{jk}^{e,f}\otimes h_{ki}^{f,f}$.  We write $\{w^1,\ldots,w^r\}$ and $\{m^1,\ldots,m^t\}$ for the dual bases for the left duals of $W$ and $M$ in $\comod_{\fd}(\mathcal{SL}_m(e))$ and in $\comod_{\fd}(\mathcal{SL}_m(f))$, respectively. It is straightforward to check that $\rho(w^i)=\sum_j w^j\otimes S_{e,e}(h_{ij}^{e,e})$ and $\rho(m^i)=\sum_j m^j\otimes S_{f,f}(h_{ij}^{f,f})$. Similarly, we have $(h_{ij}^{f,e})_{1\leq i\leq t,1\leq j\leq r}\in \mathcal{SL}_m(f,e)$ such that $\sum_j w^j\otimes S_{f,e}(h^{f,e}_{ij})\in W^*\cotens_{\mathcal{SL}_m(e)}\mathcal{SL}_m(e,f)$ forms a basis for $F(W^*)$ for $1\leq i\leq t$.

Before showing \eqref{pivS} commutes, we derive an explicit formula for $\widetilde{\xi}: F(W^*)\cong F(W)^*$ at the basis level. 
Define a $\kk$-linear map $\xi': F(W^*)\to F(W)^*$ via $\sum_{1 \leq j \leq r} w^j\otimes S_{f,e}(h^{f,e}_{ij})\mapsto m^i$ for any $1\leq i\leq t$. We claim the following diagram commutes
\[\small
\xymatrix{
F(W^*)\otimes F(W)\ar[ddd]_-{\xi}\ar[rrr]^-{\xi'\otimes \id }\ar@{=}[dr]&&&F(W)^*\otimes F(W)\ar@{=}[dl]\ar[ddd]^-{{\rm ev}}\\
&(W^*\cotens_{\mathcal{SL}_m(e)} \mathcal{SL}_m(e,f))\otimes(W\cotens_{\mathcal{SL}_m(e)} \mathcal{SL}_m(e,f))\ar[r]^-{\xi'\otimes \id}\ar[d]_-{\xi} &M^*\otimes M\ar[d]^-{{\rm ev}}&\\
&(W^*\otimes W)\cotens_{\mathcal{SL}_m(e)} \mathcal{SL}_m(e,f)\ar[r]^-{{\rm ev}\cotens_{\mathcal{SL}_m(e)}\mathcal{SL}_m(e,f)}&\kk &\\
F(W^*\otimes W)\ar@{=}[ur]\ar[rrr]^-{F({\rm ev})}&&& \kk. \ar@{=}[ul]
}
\]
In detail, the composition of the left column and bottom row maps of the inner rectangle yields that, 
\begin{align*}
    F({\rm ev})\xi \left( \left (\sum_k w^k\otimes S_{f,e}(h_{ik}^{f,e}) \right) \otimes \left (\sum_l w_l \otimes h_{lj}^{e,f} \right )\right)&\,=F({\rm ev})\left(\sum_{k,l} (w^k\otimes w_l) \otimes S_{f,e}(h_{ik}^{f,e})h_{lj}^{e,f}\right)\\
    &=\, \sum_{k,l} \langle w^k,w_l\rangle S_{f,e}(h_{ik}^{f,e})h_{lj}^{e,f}\\
    &=\, \sum_k S_{f,e}(h_{ik}^{f,e})h_{kj}^{e,f}=\varepsilon_f(h_{ij}^{f,f})=\delta_{ij}.
\end{align*}
Meanwhile, the composition of the top row and right column maps of the inner rectangle yields that
\begin{align*}
    \left\langle \xi' \left (\sum_k w^k\otimes S_{f,e}(h_{ik}^{f,e}) \right), \sum_l w_l \otimes h_{lj}^{e,f}\right\rangle=\langle m^i,m_j\rangle=\delta_{ij}.
\end{align*}
This proves our claim.  By the definition of the duality transformation \eqref{xitilde}, we observe that $\widetilde{\xi}: F(W^*)\to F(W)^*$ is the unique map that makes the above diagram commute. Hence, $\widetilde{\xi}$ and $\xi'$ coincide and we have $\widetilde{\xi}(\sum_j w^j\otimes S_{f,e}(h^{f,e}_{ij}))=m^i$ for any $1\leq i\leq t$.

Thus, tracing along the left column and bottom row of the diagram \eqref{pivS}, we have
\begin{align*}
  \left \langle \widetilde{\xi}^*j_{F(W)}(m_i), \sum_l w^l\otimes S_{f,e}(h_{jl}^{f,e})\right\rangle&\,=\left\langle \sum_k m_k\Phi_{f}^{-1}(h_{ki}^{f,f}), \widetilde{\xi}\left(\sum_l w^l\otimes S_{f,e}(h_{jl}^{f,e})\right )\right\rangle \\
  &\,= \left \langle \sum_k m_k\Phi_{f}^{-1}(h_{ki}^{f,f}), m^j \right\rangle=\sum_k\Phi_f^{-1}(h_{ki}^{f,f})\langle m_k,m^j\rangle\\
  &\,=\Phi_f^{-1}(h_{ji}^{f,f}).
\end{align*}
While tracing along the top row and right column of the diagram \eqref{pivS}, we have
\begin{align*}
&\left\langle \widetilde{\xi}F(j_W)\left (\sum_k w_k\otimes h_{ki}^{e,f}\right), \sum _lw^l\otimes S_{f,e}(h_{jl}^{f,e})\right\rangle \\ 
&\qquad = \left\langle \widetilde{\xi}\left(\sum_{k,s} w_s\otimes \Phi_e^{-1}(h_{sk}^{e,e})h_{ki}^{e,f}\right), \sum _lw^l\otimes S_{f,e}(h_{jl}^{f,e})\right\rangle\\
&\qquad = \left\langle \widetilde{\xi}\left (\sum_{s} w_s\otimes (\Phi_e^{-1}*\id)(h_{si}^{e,f})\right), \sum _lw^l\otimes S_{f,e}(h_{jl}^{f,e})\right\rangle\\
&\qquad = \left\langle \widetilde{\xi}\left (\sum_{s} w_s\otimes (S_{f,e}S_{e,f}*\Phi_f^{-1})(h_{si}^{e,f})\right ), \sum _lw^l\otimes S_{f,e}(h_{jl}^{f,e})\right\rangle\\
&\qquad = \left\langle \widetilde{\xi}\left (\sum_{k,s} w_s\otimes S_{f,e}S_{e,f}(h_{sk}^{e,f})\Phi_f^{-1}(h_{ki}^{f,f}) \right), \sum _lw^l\otimes S_{f,e}(h_{jl}^{f,e})\right\rangle\\
&\qquad = \sum_{k,s,l}\langle w_s,w^l\rangle S_{f,e}S_{e,f}(h_{sk}^{e,f}))\Phi_f^{-1}(h_{ki}^{f,f})S_{f,e}(h_{jl}^{f,e})\\
&\qquad = \sum_{k,s,l}\langle w_s,w^l\rangle S_{f,e}(h_{jl}^{f,e}S_{e,f}(h_{sk}^{e,f}))\Phi_f^{-1}(h_{ki}^{f,f})\\
&\qquad = \sum_{k}S_{f,e} \left (\sum_s h_{js}^{f,e} S_{e,f}(h_{sk}^{e,f})\right)\Phi_f^{-1}(h_{ki}^{f,f})\\
&\qquad = \sum_{k}S_{f,e}(\varepsilon_{f}(h_{jk}^{f,f}))\Phi_f^{-1}(h_{ki}^{f,f}) \\
&\qquad = \Phi_f^{-1}(h_{ji}^{f,f}).
\end{align*}
This shows that $F$ preserves the copivotal structures. As a consequence, $F$ preserves quantum dimensions, that is, $d_e(W)=d_f(F(W))$ for any $W\in \comod_{\fd}(\mathcal{SL}_m(e))$ (see \cite[Lemma 6.1]{Ng-Schauenburg}). Hence, 
\[q_{A(e,N)}(t)=\sum_{i\ge 0} d_e(A(e,N)_i)\,t^i=\sum_{i\ge 0} d_f(F(A(e,N)_i))\,t^i=\sum_{i\ge 0} d_f(A(f,N)_i)\,t^i=q_{A(f,N)}(t). \qedhere
\]
\end{proof}

\begin{remark}
 Let $\mathcal C$ be any connected copivotal cogroupoid. For any objects $X,Y\in \mathcal C$, one can similarly show that the monoidal equivalence 
 \[
 F(-)=(-)\cotens_{\mathcal C(X,X)}\mathcal C(X,Y): \comod_{\fd}(\mathcal C(X,X))\to \comod_{\fd}(\mathcal C(Y,Y))
 \]
 preserves pivotal structures, that is, the diagram
 \begin{align*}
\xymatrix{
F(W)\ar[d]_-{j_{F(W)}}\ar[rr]^-{F(j_W)} && F(W^{**})\ar[d]^-{\widetilde{\xi}}\\
F(W)^{**}\ar[rr]^-{\widetilde{\xi}^*} && F(W^*)^*
}
\end{align*}
commutes for any $W\in \comod_{\fd}(\mathcal C(X,X))$. The pivotal structures on $\comod_{\fd}(\mathcal C(X,X))$ and $\comod_{\fd}(\mathcal C(Y,Y))$ are given by their characters $\Phi_X,\Phi_Y$  as in \eqref{jW}, respectively.

As an application, in \cite{Ng-Schauenburg,NgSchauenburg2010}, higher (generalized) Frobenius--Schur indicators are defined for any $\kk$-linear $\Hom$-finite monoidal category with pivotal structure $j: \id \to (-)^{**}$. Our result together with \cite[Corollary 4.4]{Ng-Schauenburg} and (\cite[Remark 2.2]{NgSchauenburg2010}) shows that the values of those indicators in $\comod_{\fd}(\mathcal{C}(X,X))$ and $\comod_{\fd}(\mathcal{C}(Y,Y))$ are preserved under the monoidal equivalence $F$.   
\end{remark}

We wonder if the converse of \Cref{prop:SameH} holds. We provide some evidence for the case when $N=m=2$. We emphasize that $\GLQS$ equivalence generally does not imply $\SLQS$ equivalence.

\begin{proposition}
\label{equivm2}
For any nondegenerate twisted superpotentials $e\in U^{\otimes 2}$ and $f \in V^{\otimes 2}$, the following are equivalent.
    \begin{itemize}
        \item[(1)] $A(e,2)=TU_e/(\kk e)$ and $A(f,2)=TV_f/(\kk f)$ are $\SLQS$ equivalent;
        \item[(2)] $\mathcal{SL}_2(e,f)\neq 0$;
        \item[(3)] $d_e(U_e)=d_f(V_f)$;
        \item[(4)] $q_{A(e,2)}(t)=q_{A(f,2)}(t)$.
    \end{itemize}
\end{proposition}

\begin{proof}
There is a straightforward presentation of the cogroupoid $\mathcal{SL}_2$. Let $e\in U^{\otimes 2}$ and write $\mathbb E=(e_{ij})\in M_p(\kk)$ with $p=\dim U$. One can check directly that $(e,\mathbb P)$ is a nondegenerate twisted superpotential if and only if $\mathbb E$ is invertible or $\mathbb E\in \GL_p(\kk)$ with $\mathbb P=\mathbb E^T\mathbb E^{-1}$. The same holds for $f \in V^{\otimes 2}$ with $q=\dim V$. So for two nondegenerate superpotentials $e\in U^{\otimes 2}$ and $f\in V^{\otimes 2}$, write $\mathbb E\in \GL(U)$ and $\mathbb F\in \GL(V)$ for their corresponding invertible matrices. We can represent $\mathcal{SL}_2(e,f)$ with generating space $\mathbb A=(a_{ij})_{1\leq i\leq p,1\leq j\leq q}$ subject to relations 
\[
\mathbb A \mathbb F\mathbb A^T\mathbb E^{-1}=\mathbb I_{p\times p}\quad \text{and}\quad \mathbb F \mathbb A^T\mathbb E^{-1}\mathbb A=\mathbb I_{q\times q}
\]
with $\mathbb B$ implicitly given by $\mathbb F\mathbb A^T\mathbb E^{-1}$. Moreover, for any vector spaces $U,V,W$ and any nondegenerate twisted superpotentials $e\in U^{\otimes 2}$, $f\in V^{\otimes 2}$ and $g\in W^{\otimes 2}$, we have algebra maps
\[\Delta_{e,g}^f: \mathcal{SL}_2(e,f)\to \mathcal{SL}_2(f,g)\otimes \mathcal{SL}_2(g,f),\qquad  \varepsilon_e: \mathcal{SL}_2(e,e)\to \kk,\qquad  S_{e,f}: \mathcal{SL}_2(e,f)\to \mathcal{SL}_2(e,f)^{\op},\]
such that 
\begin{gather*}
\Delta_{e,f}^g(a_{ij}^{e,g})=\sum_{k=1}^q a^{e,f}_{ik}\otimes a^{f,g}_{kj},\qquad \varepsilon_{e}(a_{ij}^{e,e})=\delta_{ij},\qquad S_{e,f}\left(\mathbb A^{e,f}\right)=\mathbb E\mathbb A^{f,e}\mathbb F^{-1},
\end{gather*} 
for all possible $i,j$. Then it is straightforward to check that our cogroupoid $\mathcal {SL}_2$ is essentially the same as the cogroupoid $\mathcal B$ constructed in \cite[Definition 3.3]{Bichon2014} with $\mathcal{SL}_2(e,f)=\mathcal B(\mathbb E^{-1},\mathbb F^{-1})$. 

(1) $\Leftrightarrow$ (2): By \cite{Zhang1998}, we know $A(e,2)$ is Koszul AS-regular of dimension two as long as $\mathbb E$ is invertible. Then the equivalence follows from \Cref{thm:SLQS}.

(2) $\Leftrightarrow$ (3): Since $\mathbb P=\mathbb E^T\mathbb E^{-1}$, the character $\Phi_e: \mathcal{SL}_2(e)\to \kk$ is given by $\Phi_e(\mathbb A)=\mathbb P^T=\mathbb E^{-T} \mathbb E$ and $\Phi_e(\mathbb B)=\mathbb P^{-T}=\mathbb E^{-1}\mathbb E^T$. Thus, 
  \begin{align*}
    d_e(U_e)&\,=\sum_{i} \langle j_{U_e}(u_i),u^i\rangle= \sum_{i} \left \langle \sum_ju_j\Phi_e^{-1}(a_{ji}),u^i \right\rangle=\sum_{i,j} \langle u_j\Phi_e(S_{e,e}(a_{ji})),u^i\rangle\\
    &\,=\sum_{i,j} \langle u_j\Phi_e(b_{ji}),u^i\rangle=\sum_{i,j} P^{-T}_{ji}\langle u_j\mathbb ,u^i\rangle=\sum_{i} P^{-T}_{ii}={\rm tr}(\mathbb P^{-1})={\rm tr}(\mathbb E\mathbb E^{-T})={\rm tr}(\mathbb E^{-1}\mathbb E^T).
  \end{align*}
A consequence of \cite[Lemma 3.4]{Bichon2014} shows that $\mathcal{SL}_2(e,f)\neq 0$ if and only if ${\rm tr}(\mathbb E^{-1}\mathbb E^T)={\rm tr}(\mathbb F^{-1}\mathbb F^T)$, and hence if and only if $d_e(U_e)=d_f(V_f)$.

(3) $\Leftrightarrow$ (4): By \Cref{Le:qd}, we get $q_{A(e,2)}(t)=\frac{1}{1-d_e(U_e)t+t^2}$ and $q_{A(f,2)}(t)=\frac{1}{1-d_f(V_f)t+t^2}$. Thus, $q_{A(e,2)}(t)=q_{A(f,2)}(t)$ if and only if $d_e(U_e)=d_f(V_f)$.
\end{proof}

\section{Some examples of non-noetherian AS-regular algebras of dimension 3}\label{sec:Example}

In this section, we apply our results and provide conjectural new examples of non-noetherian AS-regular algebras of global dimension 3. We will produce these algebras as superpotential algebras with 4 generators subject to some quadratic relations. By \cite[Theorem 8.1]{ATV1990} and \cite[Theorem 0.3]{SZ1997}, such algebras are necessarily non-noetherian. By \Cref{infinite GK}, we do not assume finite Gelfand--Kirillov dimension in these examples.

\subsection{Adapting Mori--Smith results to four variables}

We first recall some notation and basic results from \cite{MS2017} that we adapt to our setting of four variables.

\begin{notation}
Throughout this section, $\kk$ is an algebraically closed field with char $\kk \neq 2,3$ and $V$ is a 4-dimensional $\kk$-vector space with basis $\{x_0, x_1, x_2, x_3\}$. Let $TV$ be the tensor algebra on $V$, $S_3$ the symmetric group on $3$ letters acting on $V^{\otimes 3}$, $\Sym^3$ the symmetric tensors in $V^{\otimes 3}$, $\alt^3(V)$ the alternating (or skew-symmetric) tensors in $V^{\otimes 3}$, and $S^3V$ the degree-$3$ component of the symmetric algebra $SV$ on $V$. The canonical homomorphism $TV \to SV$ given by $w \mapsto \overline{w}$ restricts to an isomorphism $\Sym^3 V \to S^3 V$, we denote its inverse by the symmetrization map $S^3 V \to \Sym^3 V$ defined by $f \mapsto \widehat{f}$. 
\end{notation}

The rest of this subsection is essentially adapted from \cite[Section 2.2]{MS2017} to our setting. Suppose $\kk S_3$ denotes the group algebra of $S_3$. We write \textbf{1}, \textbf{sgn}, and \textbf{2} for the trivial representation, sign representation, and 2-dimensional irreducible representation of $\kk S_3$, respectively. The elements $c$ and $s$ in $\kk S_3$ are idempotents such that 
\[c = \frac{1}{3}\big(1 + (123) + (321)\big) \qquad \text{ and } \qquad s = \frac{1}{2} c \cdot \big(1 - (12)\big) = \frac{1}{2} c \cdot \left(1 - (13)\right) = \frac{1}{2} c \cdot \left(1 - (23) \right).\] 
Moreover, we have the following: $sg = s$ for all $g \in S_3$, $cs = sc = s$, and $(1 2 3)c = c$. It follows that the elements $1-c$, $c-s$, and $s \in S_3$ form a complete set of mutually orthogonal central idempotents in $\kk S_3$. Hence, if $M$ is a left $\kk S_3$-module, then
\begin{align*}
sM &= \text{the sum of all submodules of } M \text{ that are isomorphic to \textbf{1}},\\
(c-s)M &= \text{the sum of all submodules of } M \text{ that are isomorphic to \textbf{sgn}},\\
(1-c)M &= \text{the sum of all submodules of } M \text{ that are isomorphic to \textbf{2}},\\
cM &= sM \oplus (c-s)M,\\
m &= s(m) + (c-s)(m) + (1-c)(m), \quad \text{for all } m \in M.
\end{align*}

Now for $M=V^{\otimes 3}$ where $V$ is $4$-dimensional with basis $\{x_0, x_1, x_2, x_3\}$, and for any $0 \leq i,j,k \leq 3$, we omit the tensor notation and define the map 
\[w^{(ijk)}: V^{\otimes 3} \to (c-s)V^{\otimes 3} \quad \text{ by } \quad x_ix_jx_k \mapsto \frac{1}{3}(x_ix_jx_k + x_kx_ix_j + x_jx_kx_i) - \frac{1}{3}(x_kx_jx_i + x_ix_kx_j + x_jx_ix_k),\] 
and define 
\begin{equation}
\label{w0}
   w_0 \coloneqq a_0w^{(012)} + a_1w^{(023)} + a_2w^{(013)} + a_3w^{(123)} \in V^{\otimes 3},  
\end{equation}
for some $a_0, a_1, a_2, a_3 \in \kk$. 

\begin{lemma}
\label{alt3}
Let $\{x_0, x_1, x_2, x_3\}$ be a basis for a $\kk$-vector space $V$ and retain the above notation. Then 
\[\alt^3(V) = {\rm Span}_\kk \{w^{(ijk)}\}_{0 \leq i,j,k \leq 3} = (c-s)V^{\otimes 3} \qquad \text{ and } \qquad cV^{\otimes 3} = \alt^3(V) \oplus \Sym^3(V).\]
\begin{proof}
The argument in the proof of \cite[Lemma 2.1]{MS2017} holds similarly for $\dim_{\kk} (V) = 4$.
\end{proof}
\end{lemma}

We note that the results in \cite[Section 2.3 and Section 3]{MS2017} hold as well for our setting. In particular, by \cite[Lemma 2.3]{MS2017}, take $f = v_1 \otimes v_2 \otimes v_3 \in S^3 V$ where each $v_i \in V$, its symmetrization $\widehat{f} = \frac{1}{6}\sum_{\sigma \in S_3}v_{\sigma(1)} \otimes v_{\sigma(2)} \otimes v_{\sigma(3)}$, and \cite[Theorem 3.2]{MS2017} remains true when we take $\{x_0, x_1, x_2, x_3\}$ to be the basis for $V$. We have the following description for the superpotentials in ${\rm CY}(3,2,4)$ (cf.~\eqref{CYclass}) that yield Koszul
AS-regular algebras that are Calabi--Yau of dimension $3$.

\begin{lemma}
\label{wCY}
Let $\{x_0, x_1, x_2, x_3\}$ be a basis for a $\kk$-vector space $V$ and retain the above notation. Then, every superpotential in ${\rm CY}(3,2,4)$ is given by $e= w_0 + \widehat{f}$, for some $f \in S^3 V$. In particular, $e=\widehat{f}$ if and only if $w_0=0$ if and only if $c(e)=s(e)$.
\end{lemma}

\begin{proof}
For any $e \in {\rm CY}(3,2,4)$, $e = c(e) = (c-s)e + se = w_0 + \widehat{\overline{se}}$, where the last equality follows from \Cref{alt3}, the definition of $w_0$ in \eqref{w0}, and the fact that $\overline{(-)}$ and $\widehat{(-)}$ are inverses of each other. Hence, $e=w_0 + \widehat{f}$ with $f=\overline{se} \in S^3V$. The in particular statement is straightforward. 
\end{proof}

\subsection{Examples}

In this subsection, we illustrate our results by checking the non-vanishing of the associated bi-Galois object given by two superpotentials. In particular, by \Cref{ASSuper} we know that the quantum-symmetric equivalence class of the polynomial ring is completely determined by the non-vanishing of the algebras $\mc{GL}_m(e,f)$. 

In the examples below, for particular choices of $e= w_0 + \widehat{f}$, with $w_0$ given in \eqref{w0}, defined on a 4-dimensional vector space with basis $\{x_0,x_1,x_2,x_3\}$ as described in the previous subsection, the resulting superpotential algebra $A(e,2)$ is generated by 4 generators, subject to quadratic relations. Recall that two such superpotential algebras corresponding to superpotentials $e$ and $f$ are in the same connected component of $\mc{GL}_3$ if and only if $\mc{GL}_{3}(e,f) \not \cong 0$. We write ${\rm ASreg}_3$ for the connected component of $\mc{GL}_3$ containing the polynomial ring in $3$ variables, which by \Cref{ASSuper} consists precisely of the superpotential algebras which are Koszul AS-regular of dimension $3$.

Set $f_{\text{poly}}$ as the superpotential corresponding to the polynomial ring in 3 variables, and let $e$ be an arbitrary superpotential as above. To determine whether $\mc{GL}_3(e,f_{\text{poly}})$ is nonzero, we use explicit generators and relations of the algebra $\mc{GL}_3(e,f_{\text{poly}})$ as given in \eqref{eq:alg} and attempt to compute a Groebner basis for $\mc{GL}_3(e,f_{\text{poly}})$ via Magma. In the case that the algebra is 0, the Groebner basis algorithm typically terminates very quickly; by \Cref{ASSuper}, we know in this case that $A(e,2)$ is not AS-regular. On the other hand, we find that when $\mc{GL}_3(e,f_{\text{poly}}) \not \cong 0$, the Groebner basis algorithm typically does not terminate; even when we restrict Magma to compute up to a certain degree, the number of resulting basis elements seems to grow exponentially as the degree increases. We view this case, when the algorithm does not terminate quickly, as strong evidence for the AS-regularity of $A(e,2)$, although it does not provide a rigorous proof (since a priori we do not know as a fact that the algorithm does not terminate at some later stage). For the non-AS-regular superpotential algebras, we then test them against one another to make a conjecture about whether they are in the same connected component of $\mc{GL}_3$, by the same process as above. 

By \Cref{wCY}, in the following examples, for a vector space $V$ with basis $\{x_0, x_1, x_2, x_3\}$, we compute the corresponding superpotential $e=w_0 + \lambda \widehat{\mathsf{F}}$ for some $\lambda \in \kk$ and $\mathsf{F} \in S^3V$ coming from the classification of cubic surfaces in $\mathbb P^3$. We check the (non)vanishing of $\mc{GL}_3(e,f_{\text{poly}})$ on Magma as described above and consider different values for $a_0, \ldots, a_3$ and $\lambda \in \kk$ in these examples. 

\begin{Ex}
By \cite[Theorem 2 and Table 2]{Sakamaki2010}, normal singular cubic surfaces in $\mathbb P^3$ over the field of complex numbers $\mathbb C$ are defined by $\mathsf{F} \coloneqq x_3f_2(x_0,x_1,x_2) - f_3(x_0,x_1,x_2)$ of degree $3$, where $f_2,f_3$ are given in the second and third columns of \Cref{tab:Singular cubic} below, according to the types of singularities characterized by Coxeter-Dynkin diagrams. Note that the types of singularities appearing in the first column of \Cref{tab:Singular cubic} refer to the differing number of parameters and of lines of the cubic surfaces. When $a_0=a_1=a_2=a_3=0$ (that is, when $w_0=0)$, upon rescaling of $\lambda$, the corresponding superpotential $e=\widehat{\mathsf{F}}$ is nondegenerate for all examples, with generic parameters, given in \Cref{tab:Singular cubic}. 
\begin{center}
\begin{table}[h!]
\begin{tabular}{ |p{2cm}|p{2cm}|p{7cm}|c|  }
 \hline
 Singularities& $f_2(x_0,x_1,x_2)$ &$f_3(x_0,x_1,x_2)$&Connected component\\
 \hline
 $A_1$ &  $x_0x_2-x_1^2$   & $(x_0-\alpha x_1)(-x_0+(\beta+1)x_1-\beta x_2)(x_1-\gamma x_2)$& ${\rm ASreg}_3$  \\
  $2A_1$  &  $x_0x_2-x_1^2$   & $(x_0 - 2x_1+x_2)(x_0-\alpha x_1)(x_1-\beta x_2)$&  ${\rm ASreg}_3$ \\

 $A_1 A_2$ &   $x_0x_2-x_1^2$  & $(x_0-x_1)(-x_1+x_2)(x_0-(\alpha+1)x_1+\alpha x_2)$ & ${\rm ASreg}_3$  \\

$3A_1$  &   $x_0x_2-x_1^2$  & $x_0 x_2 (x_0-(\alpha+1)x_1+\alpha x_2)$ &  ${\rm ASreg}_3$ \\

$A_1 A_3$  &  $x_0x_2-x_1^2$   & $(x_0-x_1)(-x_1+x_2)(x_0-2x_1+x_2)$ &  ${\rm ASreg}_3$ \\

 $2A_1 A_2$ &  $x_0x_2-x_1^2$   & $x_1^2(x_0-x_1)$& ${\rm ASreg}_3$  \\

$4A_1$  &  $x_0x_2-x_1^2$   & $(x_0-x_1)(x_1-x_2)x_1$ &  ${\rm ASreg}_3$ \\

$A_1 A_4$  &    $x_0x_2-x_1^2$ & $x_0^2 x_1$ & ${\rm ASreg}_3$  \\

$2A_1 A_3$  &   $x_0x_2-x_1^2$  & $x_0 x_1^2$&  ${\rm ASreg}_3$ \\

 $A_1 2A_2$ &   $x_0x_2-x_1^2$  & $x_1^3$&  ${\rm ASreg}_3$ \\

 $A_1 A_5$ &  $x_0x_2-x_1^2$   & $x_0^3$& ${\rm ASreg}_3$  \\

 $A_2$ &   $x_0x_1$  & $x_2(x_0+x_1+x_2)(\delta x_0+ \eta x_1-\delta \eta x_2)$ &  ${\rm ASreg}_3$ \\

 $2A_2$ & $x_0x_1$    & $x_2(x_1+x_2)(-x_1+\delta x_2)$ & ${\rm ASreg}_3$  \\

 $3A_2$ &  $x_0x_1$   & $x_2^3$ &  B \\

 $A_3$ &   $x_0x_1$  & $x_2(x_0+x_1+x_2)(x_0-ux_1)$&  ${\rm ASreg}_3$ \\

$A_4$  &   $x_0x_1$  & $x_0^2x_2+x_1^3-x_1 x_2^2$& ${\rm ASreg}_3$  \\
 $A_5$ &   $x_0x_1$  & $x_0^3+x_1^3-x_1 x_2^2$ & ${\rm ASreg}_3$  \\

 $D_4(1)$ &   $x_0^2$  & $x_1^3+x_2^3$& C  \\

$D_4(2)$  & $x_0^2$    & $x_1^3+x_2^3+x_0 x_1 x_2$&  D \\

$D_5$  & $x_0^2$    & $x_0 x_2^2+x_1^2x_2$ &  D \\

$E_6$  &  $x_0^2$   & $x_0 x_2^2+x_1^3$ &  C \\

 \hline
\end{tabular} 
\caption{Conjectured connected components for certain superpotential algebras corresponding to $e=\widehat{ x_3f_2 - f_3}$, related to normal singular cubic surfaces}
\label{tab:Singular cubic}
\end{table}
\end{center}
Here in \Cref{tab:Singular cubic}, $\alpha,\beta,\gamma$ are three distinct elements of $\mathbb C \setminus \{0,1\}$, $\delta, \eta$ are elements of $\mathbb C \setminus \{0,-1\}$, and $u$ is an element of $\mathbb C^{\times} \coloneqq \mathbb C \setminus \{0\}$. We pick generic parameters for our Magma computations and use letters $B,C,D$ to denote different connected components other than the AS-regular connected component ${\rm ASreg}_3$ in the last column of \Cref{tab:Singular cubic}. The connected components $B, C,$ and $D$ are all different from ${\rm ASreg}_3$. For example, this reflects the fact that $\mc{GL}_3(e,f_{\text{poly}}) \cong 0$ for $e=\widehat{\mathsf{F}}$ in the $B,C$ or $D$ connected component; whereas $\mc{GL}_3(e,f) \not \cong 0$ (conjecturally) for $e$ and $f$ labeled in the same connected component. Here, by \Cref{ASSuper} and \Cref{AS2cocyle}, we claim the following result.  
\end{Ex}

\begin{proposition}
    The superpotential algebras associated to the superpotentials $e=\widehat{ x_3f_2 - f_3}$ labeled by $B,C,D$ in the last column of  \Cref{tab:Singular cubic} are not AS-regular.
\end{proposition}

\begin{Ex}[Singularity case $A_1A_3$]
For the superpotential algebras that belong to the connected component labeled ${\rm ASreg}_3$ in the last column of \Cref{tab:Singular cubic}, we claim that they are AS-regular of dimension $3$. As an example, we show that the superpotential algebra $A(e,2)$ corresponding to $e=\widehat{ x_3f_2 - f_3}$ for the singularity case $A_1A_3$ is AS-regular. Here, $f_2=x_0x_2-x_1^2$ and $f_3=(x_0-x_1)(-x_1+x_2)(x_0-2x_1+x_2)$. The defining relations of $A(e,2) = \kk\langle x_0, x_1, x_2, x_3\rangle / (r_0, r_1, r_2, r_3)$ are:
\begin{align*}
    r_0\colon \quad & 2(x_0 x_1 + x_1 x_0) - 2(x_0 x_2 + x_2 x_0) - 6x_1^2 + 4(x_1 x_2 + x_2 x_1) - 2x_2^2 + x_2 x_3 + x_3 x_2 = 0, \\
    r_1\colon \quad & x_0^2 - 3(x_0 x_1 + x_1 x_0) + 2(x_0 x_2 + x_2 x_0) + 6x_1^2 - 3(x_1 x_2 + x_2 x_1) - (x_1 x_3 + x_3 x_1) + x_2^2 = 0, \\
    r_2\colon \quad & -2x_0^2 + 4(x_0 x_1 + x_1 x_0) - 2(x_0 x_2 + x_2 x_0) + (x_0 x_3 + x_3 x_0) - 6x_1^2 + 2(x_1 x_2 + x_2 x_1) = 0, \\
    r_3\colon \quad & x_0 x_2 + x_2 x_0 - 2x_1^2 = 0. 
\end{align*}
We observe that $A(e,2)$ is an Ore extension $R[x_3; \sigma, \delta]$ of the algebra
\[
R:=\kk \langle x_0,x_1,x_2 \rangle/(x_0 x_2 + x_2 x_0 - 2x_1^2),
\]
where $\sigma: R \to R$ is a ring homomorphism and $\delta: R \to R$ is a $\sigma$-derivation of $R$, defined as
\begin{align*}
    \sigma(x_i)&= -x_i, \quad \text{ for }i=0,1,2,\\
    \delta(x_0)&= 2x_0^2 - 4(x_0 x_1 + x_1 x_0) + 2(x_0 x_2 + x_2 x_0) + 6x_1^2 - 2(x_1 x_2 + x_2 x_1), \\
    \delta(x_1)&= x_0^2 - 3(x_0 x_1 + x_1 x_0) + 2(x_0 x_2 + x_2 x_0) + 6x_1^2 - 3(x_1 x_2 + x_2 x_1) + x_2^2, \\
    \delta(x_2)&= -2(x_0 x_1 + x_1 x_0) + 2(x_0 x_2 + x_2 x_0) + 6x_1^2 - 4(x_1 x_2 + x_2 x_1) + 2x_2^2.
\end{align*}
Note that $R$ is itself a superpotential algebra, and is AS-regular by \cite[Corollary 3.2.4]{HNUVVW2024-3}. By \cite[Theorem 0.2]{LiuWangWu2014}, an Ore extension of a twisted Calabi--Yau algebra is again twisted Calabi--Yau. Combining the result of \cite[Lemma 1.2]{RRZ2014}, this implies that $A(e,2)$ is also AS-regular.

A similar argument can be done for algebras in \Cref{tab:Singular cubic} with the same $f_2=x_0x_2-x_1^2$. One can show that they are Ore extensions $R[x_3; \sigma, \delta]$ of the same $R$ and $\sigma$ as above with other suitable choices of $\delta$, showing they are AS-regular. 
\end{Ex}

\begin{Ex}
When $w_0=0$, the superpotential  $e=\widehat{\mathsf{F}}$ is degenerate for the singularity case $\widetilde{E}_6$ from \cite[Table 2]{Sakamaki2010}, hence we do not include it in \Cref{tab:Singular cubic}.
\end{Ex}

\begin{Ex}
\label{ex-wzero}
    In the case when $e=w_0+ \lambda \widehat{\mathsf{F}} = w_0+ \lambda \widehat{(x_3f_2 - f_3)}$, where $w_0$ is given in \eqref{w0}, for generic scalars $a_0, a_1,a_2,a_3,\lambda \in \mathbb C^{\times}$, by the Magma computation we conjecture that $\mc{GL}_3(e,f_{\text{poly}}) \not \cong 0$ and all superpotential algebras related to singular cubic surfaces given in \cite[Table 2]{Sakamaki2010}, including the singularity case $\widetilde{E}_6$, belong to the AS-regular connected component ${\rm ASreg}_3$.  
    
    However, when at least one of the scalars $\{a_0, a_1, a_2, a_3\}$ is zero, we do not necessarily expect to have similar (non)degeneracy of $e=w_0+ \lambda \widehat{\mathsf{F}}$ and non-vanishing of $\mc{GL}_3(e,f_{\text{poly}})$ in these cases. For example, for the superpotential $e=w_0+ \lambda \widehat{\mathsf{F}}$ corresponding to the singularity case $D_5$ on the second to last row of Table 1, when $a_0=1, a_1=a_2=a_3=0$ and $\lambda=1$, our computation gives $\mc{GL}_3(e,f_{\text{poly}})=0$; whereas when $a_1=1, a_0=a_2=a_3=0$ and $\lambda=1$, Magma does not terminate and we conjecture that $\mc{GL}_3(e,f_{\text{poly}}) \not \cong 0$. 
\end{Ex}

\begin{Ex}
    When $\lambda = 0$, the superpotential $e$ corresponding to the polynomial $w_0$ is degenerate, and hence the superpotential algebra $A(e,2)$ cannot be studied using the cogroupoid $\mc{GL}_3$. 
\end{Ex}

\begin{Ex}
        We are not restricted to only considering superpotentials arising from singular surfaces. If we take the superpotential arising from the smooth Clebsch surface with defining polynomial
    \[
    \mathsf{F}=x_0^3+x_1^3+x_2^3+x_3^3 - (x_0+x_1+x_2+x_3)^3
    \]
    and consider $e=\widehat{\mathsf{F}}$ (that is, when $w_0=0$), we obtain a superpotential algebra which appears to be in the AS-regular connected component ${\rm ASreg}_3$. On the other hand, if we take a superpotential arising from the smooth Fermat surface with defining polynomial 
    \[
    \mathsf{F}=x_0^3+x_1^3+x_2^3+x_3^3
    \]
    and consider $e=\widehat{\mathsf{F}}$, we get a superpotential algebra which appears to be in the $C$ connected component. In this case, we see from \Cref{tab:Singular cubic} that singular cubic surfaces of types $D_4(1)$ and $E_6$ are in the same connected component $C$ as the Fermat smooth cubic surface. 
\end{Ex}


\bibliography{Superpotential}
\bibliographystyle{myamsalpha}
    
\end{document}